\documentclass[12pt,a4paper]{article}
\usepackage{amssymb}
\usepackage{amsmath}
\usepackage{graphicx}
\usepackage{hyperref}

\input{tcilatex}
\begin{document}

\title{Stepwise global error control in Euler's method using the DP853
triple and the Taylor remainder term}
\author{J. S. C. Prentice \\
Senior Research Officer\\
Mathsophical Ltd.\\
Johannesburg, South Africa}
\maketitle

\begin{abstract}
We report on a novel algorithm for controlling global error in a
step-by-step (stepwise) sense, in the numerical solution of a scalar,
autonomous, nonstiff or weakly stiff problem. The algorithm exploits the
remainder term of a Taylor expansion of the solution. It requires the use of
the DP853 triple to solve an auxiliary problem which, in turn, enables the
remainder term to be determined. A \textit{quenching} process then allows
the solution generated by Euler's method to be controlled. We have achieved
tolerances on the relative global error as strict as $10^{-10}.$
\end{abstract}



\section{Introduction}

We seek to solve numerically the real-valued scalar autonomous problem%
\begin{align}
y^{\prime }\left( x\right) & =f\left( y\right)  \label{autonomuos DE} \\
y\left( x_{0}\right) & =y_{0}  \notag \\
x\in & \left[ x_{0},x_{N}\right] ,  \notag
\end{align}%
using an algorithm based on Runge-Kutta methods, with \textit{stepwise}
(step-by-step) global error control. Work has been done on stepwise
estimation and suppression of the global error in these types of problems,
and we will discuss this in a later section. There, we will also refer to
our own earlier work in this regard. In this brief introduction, it suffices
to state that we intend to describe a novel approach to controlling the
global error in Euler's method applied to (\ref{autonomuos DE}), in a
stepwise sense, by exploiting the remainder term of a low-order Taylor
expansion. Our algorithm will make use of the DP853 embedded triple \cite%
{H+N+W}. Remarkably enough, we have been able to achieve relative global
error tolerances as strict as $10^{-10}$.

\section{Relevant Concepts}

\subsection{The Taylor-Lagrange function}

Let $y\left( x\right) $ be a real-valued univariate function, and assume
that $y\left( x\right) $ is as differentiable as is required in this paper.
Taylor's theorem \cite{Spivak} provides the following result:%
\begin{equation}
y\left( x\right) =y\left( x_{0}\right) +y^{\prime }\left( \xi _{x},\mu
_{x}\right) \left( x-x_{0}\right) =y_{0}+f\left( \mu _{x}\right) \left(
x-x_{0}\right) ,  \label{Taylor}
\end{equation}%
where $\mu _{x}=y\left( \xi _{x}\right) $ and $x_{0}<\xi _{x}<x.$ The second
term on the RHS is the \textit{remainder term}, presented here in \textit{%
Lagrange} form \cite{Apost},\cite{boman}. By differentiating (\ref{Taylor})
with respect to $x,$we find%
\begin{align}
y^{\prime }\left( x\right) & =f\left( \mu _{x}\right) +\left( x-x_{0}\right) 
\frac{\partial f\left( \mu _{x}\right) }{\partial \mu _{x}}\frac{d\mu _{x}}{%
dx}  \notag \\
& =f\left( \mu _{x}\right) +f_{y}\left( \mu _{x}\right) \left(
x-x_{0}\right) \frac{d\mu _{x}}{dx}.  \label{dy-dyx0...}
\end{align}%
Now, the LHS of this expression exists (see (\ref{autonomuos DE})), and so
we will assume that the RHS also exists. Hence, we assume that $\frac{d\mu
_{x}}{dx}$ exists. This suggests that $\mu _{x}$ is a function of $x,$ and
we will sometimes write $\mu \left( x\right) $ or simply $\mu $ in place of $%
\mu _{x}.$ We will refer to $\mu _{x}$ as the \textit{Taylor-Lagrange
function,} or simply \textit{Lagrange function, }for this particular
situation. In (\ref{Taylor}) and (\ref{dy-dyx0...}), we have used the
notation%
\begin{align*}
y^{\prime }\left( \xi _{x},\mu _{x}\right) & \equiv y^{\prime }\left( \mu
\left( x\right) \right) =f\left( \mu \left( x\right) \right) \equiv f\left(
\mu _{x}\right) \\
y^{\prime \prime }\left( \xi _{x},\mu _{x}\right) & \equiv \frac{dy^{\prime
}\left( \mu \left( x\right) \right) }{d\mu \left( x\right) }=\left. \frac{%
dy^{\prime }\left( y\right) }{dy}\right\vert _{y=\mu _{x}}=\left. \frac{df}{%
dy}\right\vert _{y=\mu _{x}}\equiv f_{y}\left( \mu _{x}\right) .
\end{align*}%
Rearranging (\ref{dy-dyx0...}) and using (\ref{autonomuos DE}) and (\ref%
{Taylor}) gives%
\begin{align}
\frac{d\mu }{dx}& =\frac{f\left( y\right) -f\left( \mu \right) }{f_{y}\left(
\mu \right) \left( x-x_{0}\right) }  \notag \\
& =\frac{f\left( y_{0}+f\left( \mu \right) \left( x-x_{0}\right) \right)
-f\left( \mu \right) }{f_{y}\left( \mu \right) \left( x-x_{0}\right) }
\label{dmu/dx=...} \\
& \equiv g\left( x,\mu \right) .  \notag
\end{align}%
This is an initial-value problem that can, in principle, be solved to yield
the Lagrange function $\mu \left( x\right) $ for a suitable initial value.
Once $\mu \left( x\right) $ is known, the remainder term $f\left( \mu
\right) \left( x-x_{0}\right) $ is easily computed. Then, adding the
remainder term to $y\left( x_{0}\right) $ yields, from (\ref{Taylor}), an
approximation to $y\left( x\right) .$ For interest's sake, we refer the
reader to our recent work in this regard \cite{Prentice Taylor}.

\subsection{Runge-Kutta Methods}

The topic of Runge-Kutta (RK) methods is very broad, indeed. We will assume
the reader has suitable familiarity with these methods, and list Butcher 
\cite{Butcher} as a good general reference. Here, we will simply establish
notation and terminology relevant to this paper.

To solve $\mu ^{\prime }=g\left( x,\mu \right) $ numerically, we partition $%
\left[ x_{0},x_{N}\right] $ using the nodes $x_{0}<x_{1}<\ldots <x_{N}$ and
we use the RK method%
\begin{equation*}
\mu _{i+1}=\mu _{i}+h_{i+1}\Phi \left( x_{i},\mu _{i}\right) .
\end{equation*}%
Here, $\mu _{i}$ and $\mu _{i+1}$ are the approximate numerical solutions at 
$x_{i}$ and $x_{i+1},$ respectively, $h_{i+1}\equiv x_{i+1}-x_{i}$ is the 
\textit{stepsize}, and $\Phi \left( x,\mu \right) $ is known as the \textit{%
increment function}. The increment function defines the RK method. We will
make use of several RK methods in this work - denoted RK1, RK4, RK7 and
DP853 - and we will describe them when appropriate. Also, the consistency
property of RK methods means that%
\begin{align}
\lim_{h_{i+1}\rightarrow 0}\Phi \left( x_{i},\mu _{i}\right) & =g\left(
x_{i},\mu _{i}\right)  \notag \\
\lim_{h_{i+1}\rightarrow 0}\Phi _{\mu }\left( x_{i},\mu _{i}\right) &
=g_{\mu }\left( x_{i},\mu _{i}\right)  \label{consistency} \\
\lim_{h_{i+1}\rightarrow 0}\Phi _{\mu \mu }\left( x_{i},\mu _{i}\right) &
=g_{\mu \mu }\left( x_{i},\mu _{i}\right) ,  \notag
\end{align}%
where the subscript $\mu $ denotes partial differentiation with respect to $%
\mu .$

\subsubsection{The DP853 embedded triple}

The DP853 embedded triple \cite{H+N+W} comprises three RK methods (of orders
3, 5 and 8) that use differing combinations of the same internal stages to
generate their output. There are 12 stages in DP853, each of which has the
general form%
\begin{equation*}
k_{q}=h_{i+1}g\left( x_{i}+c_{q}h_{i+1},\mu
_{i}^{V}+a_{q,1}k_{1}+a_{q,2}k_{2}+\ldots +a_{q,q-1}k_{q-1}\right)
\end{equation*}%
for each $q\in \left[ 1,2,3,\ldots ,12\right] .$ The methods are explicit,
so that each stage is dependent on the previous stages, and they must be
computed sequentially. We use the superscripts $L,H$ and $V$ to denote the
third-, fifth- and eighth-order solutions, respectively. Hence, the outputs
of DP853 are given by%
\begin{eqnarray}
\mu _{i+1}^{L} &=&\mu _{i}^{L}+\sum\limits_{q=1}^{12}b_{q}^{L}k_{q}\equiv
\mu _{i}^{L}+h_{i+1}\Phi ^{L}\left( x_{i},\mu _{i}^{L}\right)
\label{mu L = mu L + ...} \\
\mu _{i+1}^{H} &=&\mu _{i}^{H}+\sum\limits_{q=1}^{12}b_{q}^{H}k_{q}\equiv
\mu _{i}^{H}+h_{i+1}\Phi ^{H}\left( x_{i},\mu _{i}^{H}\right)  \notag \\
\mu _{i+1}^{V} &=&\mu _{i}^{V}+\sum\limits_{q=1}^{12}b_{q}^{V}k_{q}.\equiv
\mu _{i}^{V}+h_{i+1}\Phi ^{V}\left( x_{i},\mu _{i}^{V}\right)  \notag
\end{eqnarray}%
wherein we have implicitly defined coefficients $b_{q}^{\cdots }$ and
appropriate increment functions. It transpires that DP853 has $b_{q}^{H,V}=0$
for $q\in \left[ 2,3,4,5\right] $ and $b_{q}^{L}=0$ for $q\in \left[
2,3,4,5,6,7,8,10,11\right] .$

We modify (\ref{mu L = mu L + ...}) slightly: we use $\mu _{i}^{H}$ as input
for $\mu _{i+1}^{L},$ i.e.%
\begin{equation}
\mu _{i+1}^{L}=\mu _{i}^{H}+\sum\limits_{q=1}^{12}b_{q}^{L}k_{q}.
\label{mu L = mu H +...}
\end{equation}%
We discuss the value of this approach in the next subsection.

\subsection{Local and Global Error}

The RK \textit{local error} at $x_{i+1}$ is%
\begin{equation*}
\delta _{i+1}\equiv \mu \left( x_{i+1}\right) -\left( \mu \left(
x_{i}\right) +h_{i+1}\Phi \left( x_{i},\mu \left( x_{i}\right) \right)
\right)
\end{equation*}%
and the RK global error at $x_{i+1}$ is%
\begin{equation*}
\Delta _{i+1}\equiv \mu \left( x_{i+1}\right) -\mu _{i+1}=\mu \left(
x_{i+1}\right) -\left( \mu _{i}+h_{i+1}\Phi \left( x_{i},\mu _{i}\right)
\right) ,
\end{equation*}%
where $\mu \left( x_{i}\right) $ and $\mu \left( x_{i+1}\right) $ are the
exact solutions at $x_{i}$ and $x_{i+1},$ respectively. The difference
between the two concepts is subtle: the local error definition requires that 
$\mu \left( x_{i}\right) $ is known.

An RK method of \textit{order} $p,$ written RK$p,$ has%
\begin{align}
\delta _{i+1}& \propto h_{i+1}^{p+1}+\ldots =L_{i+1}h_{i+1}^{p+1}+\ldots
\label{local error} \\
\Delta _{i+1}& \propto h^{p}+\ldots =G_{i+1}h^{p}+\ldots
\label{global error}
\end{align}%
In the latter expression, the symbol $h$ represents the mean stepsize over
the subinterval $\left[ x_{0},x_{i+1}\right] .$ The ellipses represent
higher-order terms.

The local and global errors are related. We have shown \cite{Prentice 1}\cite%
{Prentice 2} that%
\begin{equation}
\Delta _{i+1}=\delta _{i+1}+\Delta _{i}+h_{i+1}\Phi _{\mu }\left( x_{i},\mu
_{i}\right) \Delta _{i}+\ldots  \label{Delta i+1 = ...}
\end{equation}%
where $\Delta _{i}$ is the global error in the input $\mu _{i}.$ Hence,
ignoring higher order terms, (\ref{mu L = mu H +...}) yields%
\begin{equation*}
\Delta _{i+1}^{L}=\delta _{i+1}^{L}+\left( 1+h_{i+1}\Phi _{\mu }\left(
x_{i},\mu _{i}^{H}\right) \right) \Delta _{i}^{H}.
\end{equation*}%
Compare this with 
\begin{equation*}
\Delta _{i+1}^{L}=\delta _{i+1}^{L}+\left( 1+h_{i+1}\Phi _{\mu }\left(
x_{i},\mu _{i}^{L}\right) \right) \Delta _{i}^{L}
\end{equation*}%
which follows from (\ref{mu L = mu L + ...}). The second term on the RHS is
the contribution to $\Delta _{i+1}^{L}$ due to the global error of the input 
$\mu _{i}^{H}$ or $\mu _{i}^{L}.$ Since we expect $\left\vert \Delta
_{i}^{H}\right\vert \ll \left\vert \Delta _{i}^{L}\right\vert ,$ we see that
(\ref{mu L = mu H +...}) reduces the magnitude of $\Delta _{i+1}^{L}$
compared with what it would otherwise have been.

\subsection{Error Control}

\subsubsection{Local error}

It is possible to adjust the stepsize so as to control the local error. If $%
\varepsilon _{\delta }>0$ denotes a user-defined tolerance on the local
error, then 
\begin{equation*}
\left\vert \delta _{i+1}\right\vert =\left\vert
L_{i+1}h_{i}^{p+1}\right\vert \leqslant \varepsilon _{\delta }\Rightarrow
h_{i+1}\leqslant \left( \frac{\varepsilon _{\delta }}{\left\vert
L_{i+1}\right\vert }\right) ^{\frac{1}{p+1}}
\end{equation*}%
provides an estimate of the stepsize necessary to bound $\left\vert \delta
_{i+1}\right\vert .$ In this expression we have ignored higher-order terms
present in the expansion of $\delta _{i+1}$ (see (\ref{local error})). To
cater for this omission, it is common practice to use a \textit{safety factor%
} $0<\eta <1,$ so that we have%
\begin{equation*}
h_{i+1}=\eta \left( \frac{\varepsilon _{\delta }}{\left\vert
L_{i+1}\right\vert }\right) ^{\frac{1}{p+1}}
\end{equation*}%
as the desired value of $h_{i+1}.$ We will use $\eta =0.85$ throughout this
study. Furthermore, it is possible to control the \textit{local error
density }$\rho _{i+1},$ or \textit{local error per unit step} as it is often
known, 
\begin{equation}
\rho _{i+1}\equiv \frac{\left\vert \delta _{i+1}\right\vert }{h_{i}}%
=\left\vert L_{i+1}h_{i+1}^{p}\right\vert \leqslant \varepsilon _{\rho
}\Rightarrow h_{i+1}=\eta \left( \frac{\varepsilon _{\rho }}{\left\vert
L_{i+1}\right\vert }\right) ^{\frac{1}{p}},  \label{absolute error per step}
\end{equation}%
where $\varepsilon _{\rho }=\varepsilon _{\delta }/h_{i+1}.$

There\ are two virtues in controlling the error density. Firstly, the
exponent is larger $\left( 1>\frac{1}{p}>\frac{1}{p+1}\right) $ and,
anticipating that $\frac{\varepsilon _{\rho }}{\left\vert L_{i+1}\right\vert 
}\ll 1,$ we will find that $h_{i+1}$ will be smaller, meaning that $%
\left\vert \delta _{i+1}\right\vert \leqslant \varepsilon _{\delta }$ is
satisfied, as well.

Secondly, if $h_{i+1}$ is suitably small we find, from (\ref{Delta i+1 = ...}%
), the recursion%
\begin{equation}
\Delta _{i+1}=\delta _{i+1}+\Delta _{i}\Rightarrow \Delta _{i+1}\leqslant
\left\vert \delta _{i+1}\right\vert +\Delta _{i}.  \label{recursion}
\end{equation}%
With $\Delta _{0}=0$ (because the initial value is known exactly) and
assuming $\left\vert \delta _{i+1}\right\vert \leqslant \varepsilon _{\rho
}h_{i+1},$ we find%
\begin{equation*}
\Delta _{N}\leqslant \varepsilon _{\rho
}\sum\limits_{i=0}^{N-1}h_{i+1}=\varepsilon _{\rho }\left(
x_{N}-x_{0}\right) .
\end{equation*}%
The recursion (\ref{recursion}) is monotonically increasing, so that $%
\varepsilon _{\rho }\left( x_{N}-x_{0}\right) $ is a rough bound on the
global error that can be made \textit{a priori}. Hence, for a desired value
of $\Delta _{N},$ we could estimate the necessary local tolerance $%
\varepsilon _{\rho }$ via%
\begin{equation}
\varepsilon _{\rho }\leqslant \frac{\Delta _{N}}{x_{N}-x_{0}}.
\label{eps rho <= Delta/interval}
\end{equation}

We prefer to control \textit{relative} error, as in 
\begin{equation}
\frac{\rho _{i+1}}{\left\vert \mu _{i+1}\right\vert }\equiv \frac{\left\vert
\delta _{i+1}\right\vert }{\left\vert \mu _{i+1}\right\vert h_{i+1}}%
\leqslant \varepsilon _{\rho }\Rightarrow h_{i+1}=\eta \left( \frac{%
\varepsilon _{\rho }\left\vert \mu _{i+1}\right\vert }{\left\vert
L_{i}\right\vert }\right) ^{\frac{1}{p}},  \label{rho / mu = eta...}
\end{equation}%
when $\left\vert \mu _{i+1}\right\vert >1,$ and \textit{absolute} error, as
in (\ref{absolute error per step}), when $\left\vert \mu _{i+1}\right\vert
\leqslant 1.$ The difference between the two cases is in the numerator on
the RHS of (\ref{rho / mu = eta...}).

For relative error control, we have $\rho _{i+1}\leqslant \varepsilon _{\rho
}\left\vert \mu _{i+1}\right\vert ,$ and for absolute error control we have $%
\rho _{i+1}\leqslant \varepsilon _{\rho }.$ These two cases can be
conveniently expressed as 
\begin{equation*}
\rho _{i+1}\leqslant \varepsilon _{\rho }\max \left( 1,\left\vert \mu
_{i+1}\right\vert \right) .
\end{equation*}%
The RHS is a \textit{hinge} function (of $\left\vert \mu _{i+1}\right\vert
), $ and is somewhat different to the often-used linear expression%
\begin{equation*}
\rho _{i+1}\leqslant \varepsilon _{\rho }\left( 1+\left\vert \mu
_{i+1}\right\vert \right) .
\end{equation*}%
We see that the linear expression gives $\rho _{i+1}\leqslant \varepsilon
_{\rho }$ only when $\left\vert \mu _{i+1}\right\vert =0.$ For all other
values of $\left\vert \mu _{i+1}\right\vert ,$ the bound imposed on $\rho
_{i+1}$ is larger than $\varepsilon _{\rho }$ or, since $\rho
_{i+1}/\left\vert \mu _{i+1}\right\vert \leqslant \varepsilon _{\rho }\left(
1/\left\vert \mu _{i+1}\right\vert +1\right) ,$ the bound imposed on $\rho
_{i+1}/\left\vert \mu _{i+1}\right\vert $ is larger than $\varepsilon _{\rho
}.$ Hence, we prefer to use the hinge function, since it respects the
user-defined bound $\varepsilon _{\rho }.$

\subsubsection{Global error}

Global error control is a more difficult task, particularly in a stepwise
sense. The easiest implementable form of global error control requires an
estimation of the error after the integration has been done, and then
solving the problem again using a new, smaller stepsize. This is known as 
\textit{reintegration}, and is an \textit{a posteriori} form of control.
Much work has been done regarding global error estimation (see, for example, 
\cite{aid}\cite{calvo 1}\cite{Calvo 2}\cite{dp 1}\cite{dp 2}\cite{kulikov 1}%
\cite{kulikov 2}\cite{kulikov 3}\cite{lang}\cite{macd}\cite{shampine}\cite%
{skeel}\cite{W+K}). However, there are cases to be made for controlling the
global error as the integration proceeds \cite{Prentice 4}\cite{Prentice 7}%
\cite{Prentice 8}\cite{Prentice 9}\cite{Prentice 10}\cite{Prentice 11}\cite%
{Prentice 12}. In this regard, much less work has been done. The
quasi-consistent peer methods (see \cite{kulikov 4}\cite{kulikov 5}\cite%
{kulikov 6}\cite{weiner}\cite{W+K}\ and many references therein) have shown
very good performance, particularly when dealing with strongly stiff
systems, but these methods are not quite \textit{stepwise} control methods.
They typically require a "reboot" of the method when a certain tolerance has
been breached, which results in the method having to backtrack several,
perhaps many, nodes, resulting in what is effectively a partial
reintegration. In this paper, and in our other work \cite{Prentice 3}\cite%
{Prentice 5}\cite{Prentice 6}\cite{Prentice 7}\cite{Prentice 8}, we insist
on a backtrack of no more than a single node, when required, using a process
we refer to as \textit{quenching}.

\section{The Algorithm}

We intend to construct an algorithm that will allow for local control of the
global error in a RK1 method, by determining the remainder term in (\ref%
{Taylor}). This is a more subtle approach than the direct extrapolation
technique we have previously used \cite{Prentice 3}\cite{Prentice 5}\cite%
{Prentice 6}. We believe the best way to describe our algorithm is in a
step-by-step manner that follows the sequential implementation of the
algorithm itself.

\begin{enumerate}
\item Using $f\left( y\right) $ and (\ref{autonomuos DE}), we construct $%
g\left( x,\mu \right) .$

\item We choose a point $x_{\mu }$ very close to $x_{0},$ and use a
high-order RK method - such as RK4 or even RK7 - \ to find a solution $%
y_{\mu }$ of (\ref{autonomuos DE}) at $x_{\mu }.$ We assume $y_{\mu }$ to be
sufficiently accurate so as to be regarded as practically exact. From (\ref%
{Taylor}), we then have%
\begin{equation*}
y_{\mu }-y_{0}-f\left( \mu \right) \left( x_{\mu }-x_{0}\right) =0,
\end{equation*}%
which can be solved using a nonlinear solver to find $\mu $ (see the
Appendix; there is an important subtlety to consider). We label this value $%
\mu _{1}$ and we relabel $x_{\mu }$ as $x_{1}.$ Moreover, we observe that (%
\ref{Taylor}) is trivially satisfied when $x=x_{0},$ so that $\mu $ is
arbitrary at $x_{0}.$ There is no harm, then, in setting $\mu _{0}=0.$ This
labelling is merely a logistical convenience, and it provides the ordered
pairs $\left\{ \left( x_{0},\mu _{0}\right) ,\left( x_{1},\mu _{1}\right)
\right\} $ as the first two values in our numerical solution for $\mu .$ 
\textbf{NB:} We also relabel $y_{\mu }$ as $y_{1}.$

\item We now use a pair of RK methods to determine $\mu $ at $%
x_{2}=x_{1}+h_{2},$ with $\left( x_{1},\mu _{1}\right) $ as input. The
stepsize $h_{2}$ is user-defined (more on this point later). These RK
methods are the third- and fifth-order components of the DP853 embedded
triple. We need two values of $\mu $ at $x_{2},$ obtained with methods of
differing order, to be able to implement local error control. We will denote
these two solutions $\mu _{2}^{L}$ and $\mu _{2}^{H}$. Using DP853 sets the
order parameter $p=3$ in the steps that follow. However, other users may
choose to use RK methods of different orders, so we will persist in using
the symbol $p$ for the sake of generality.

\item The next step is to implement \textbf{local error control} (LEC) in $%
\mu $ at $x_{2}.$ We assume that $\mu _{2}^{H}$ is sufficiently accurate,
relative to $\mu _{2}^{L},$ such that%
\begin{equation*}
\mu _{2}^{H}-\mu _{2}^{L}=L_{2}h_{2}^{p+1}\Rightarrow L_{2}=\frac{\mu
_{2}^{H}-\mu _{2}^{L}}{h_{2}^{p+1}}
\end{equation*}%
is feasible. We now set $\mu _{2}=\mu _{2}^{H}$ and apply the process%
\begin{align*}
\left\vert \mu _{2}\right\vert & >1\text{ and }\left\vert
L_{2}h_{2}^{p}\right\vert >\varepsilon _{\rho }\left\vert \mu
_{2}\right\vert \Rightarrow h_{2}=\eta \left( \frac{\varepsilon _{\rho
}\left\vert \mu _{2}\right\vert }{\left\vert L_{2}\right\vert }\right) ^{%
\frac{1}{p}} \\
\left\vert \mu _{2}\right\vert & \leqslant 1\text{ and }\left\vert
L_{2}h_{2}^{p}\right\vert >\varepsilon _{\rho }\Rightarrow h_{2}=\eta \left( 
\frac{\varepsilon _{\rho }}{\left\vert L_{2}\right\vert }\right) ^{\frac{1}{p%
}}
\end{align*}%
to compute a new stepsize, if necessary. We refer to this error control as 
\textit{primary} LEC.

There is a potentially valuable refinement that we can implement at this
stage. Recall that we intend to determine the remainder term $f\left( \mu
\right) \left( x-x_{0}\right) $. If we assume $\mu _{2}\approx \mu \left(
x_{2}\right) -L_{2}h_{2}^{p+1}$ then we have%
\begin{align*}
f\left( \mu _{2}\right) \left( x_{2}-x_{0}\right) & =f\left( \mu \left(
x_{2}\right) -L_{2}h_{2}^{p+1}\right) \left( x_{2}-x_{0}\right) \\
& =f\left( \mu \left( x_{2}\right) \right) \left( x_{2}-x_{0}\right)
-L_{2}h_{2}^{p+1}f_{y}\left( \mu \left( x_{2}\right) \right) \left(
x_{2}-x_{0}\right) .
\end{align*}%
The term $L_{2}h_{2}^{p+1}f_{y}\left( \mu \left( x_{2}\right) \right) \left(
x_{2}-x_{0}\right) \approx L_{2}h_{2}^{p+1}f_{y}\left( \mu _{2}\right)
\left( x_{2}-x_{0}\right) $ can be thought of as a local error in $%
y_{0}+f\left( \mu _{2}\right) \left( x_{2}-x_{0}\right) .$ We can apply the
same reasoning as above to compute a new stepsize for a given tolerance on
this particular local error, as in%
\begin{align}
& \frac{\left\vert L_{2}f_{y}\left( \mu _{2}\right) \left( x-x_{0}\right)
\right\vert h_{2}^{p+1}}{h_{2}\max (1,\left\vert y_{0}+f\left( \mu
_{2}\right) \left( x_{2}-x_{0}\right) \right\vert )}\leqslant \varepsilon
_{\rho }  \notag \\
\Rightarrow h_{2}& =\eta \left( \frac{\varepsilon _{\rho }\max (1,\left\vert
y_{0}+f\left( \mu _{2}\right) \left( x_{2}-x_{0}\right) \right\vert )}{%
\left\vert L_{2}f_{y}\left( \mu _{2}\right) \left( x_{2}-x_{0}\right)
\right\vert }\right) ^{\frac{1}{p}}.  \label{h2 = eta...}
\end{align}%
Note the hinge function in the numerator. The division by $h_{2}$
corresponds to \textit{error per unit step} control. If this value of $h_{2}$
is smaller than the value obtained from the earlier process, then it is
chosen as the new stepsize. We refer to this error control as \textit{%
secondary} LEC. Now, given that we propagate $\mu _{1}^{H}$ in DP3 (see (\ref%
{mu L = mu H +...})) and that, for strict tolerances, stepsizes may
generally be small, $L_{2}h_{2}^{p+1}$ may be a good approximation to the
global error $\Delta _{2}.$ Consequently, (\ref{h2 = eta...}) constitutes a
form of global error control, or at least global error suppression, in $%
y_{0}+f\left( \mu _{2}\right) \left( x_{2}-x_{0}\right) .$ This is because
the quantity $\Delta _{2}f_{y}\left( \mu \left( x_{2}\right) \right) \left(
x_{2}-x_{0}\right) $ \textit{is} the global error in $y_{0}+f\left( \mu
_{2}\right) \left( x_{2}-x_{0}\right) .$ The nett result is that secondary
LEC will likely lead to small global errors in $y_{0}+f\left( \mu
_{2}\right) \left( x_{2}-x_{0}\right) ,$ a favourable and useful outcome.
Secondary LEC should be used with some caution, however. In general, the
denominator contains the factor $\left( x_{i}-x_{0}\right) $ which could
become large for large intervals. This might lead to excessively small
values for $h_{i}.$ Perhaps a lower limit on the stepsize should be
considered to counteract this. Perhaps the secondary LEC stepsize should not
be allowed to be smaller than one-tenth of the primary LEC stepsize, say. We
did not need to use such a lower limit in our numerical calculations, but it
is something that is worth considering.

Effectively, what we have done here is to attempt to control the relative
local error per step in both $\mu _{2}$ and $y_{0}+f\left( \mu _{2}\right)
\left( x_{2}-x_{0}\right) .$ We then compute a new value for $\mu _{2}$
using the new stepsize (i.e. we repeat Step 3 using the new stepsize), and
then use the new values for $\mu _{2}^{L}$ and $\mu _{2}^{H}$ so obtained to
find a new value for $L_{2}.$ In the event that no stepsize adjustment is
necessary, we simply proceed directly to step 5 (note: in the event of
repeating Step 3, we do \textit{not} then repeat Step 4, as well - we move
directly to step 5).

\item The global error in $\mu _{2}$ can now be estimated using 
\begin{equation*}
\Delta _{2}=\mu _{2}^{V}-\mu _{2}^{L}
\end{equation*}%
where $\mu _{2}^{V}$ is the eighth-order solution from DP853. We prefer to
overestimate $\Delta _{2},$ so we do not use $\Delta _{2}=\mu _{2}^{V}-\mu
_{2}^{H},$ even though we have set $\mu _{2}=\mu _{2}^{H}.$ Of course, it is
the relative global error%
\begin{equation*}
\Delta _{2}^{rel}\equiv \frac{\Delta _{2}}{\max \left( 1,\left\vert \mu
_{2}\right\vert \right) }
\end{equation*}%
that actually interests us.

\item We use Euler's Method (RK1) to find $y_{2},$ as in%
\begin{equation*}
y_{2}=y_{1}+h_{2}f\left( y_{1}\right) .
\end{equation*}%
It is, of course, the global error in $y_{2}$ that we ultimately wish to
control.

\item We find the remainder term%
\begin{equation*}
R_{2}=f\left( \mu _{2}\right) \left( x_{2}-x_{0}\right) .
\end{equation*}

\item We estimate the relative global error in the remainder term $\Delta
R_{2}$ by%
\begin{align*}
f\left( \mu _{2}\right) \left( x_{2}-x_{0}\right) & =f\left( \mu \left(
x_{2}\right) -\Delta _{2}\right) \left( x_{2}-x_{0}\right) \\
& \approx f\left( \mu \left( x_{2}\right) \right) \left( x_{2}-x_{0}\right)
-f_{y}\left( \mu _{2}\right) \left( x_{2}-x_{0}\right) \Delta _{2} \\
& \text{ \ \ \ }+\frac{f_{yy}\left( \mu _{2}\right) \left(
x_{2}-x_{0}\right) \Delta _{2}^{2}}{2} \\
\Rightarrow \Delta R_{2}& =\frac{f_{yy}\left( \mu _{2}\right) \left(
x_{2}-x_{0}\right) \Delta _{2}^{2}-2f_{y}\left( \mu _{2}\right) \left(
x_{2}-x_{0}\right) \Delta _{2}}{2\max (1,\left\vert R_{2}\right\vert )}.
\end{align*}

\item We now have the Euler solution $y_{2}$ and the Taylor approximation%
\begin{eqnarray*}
y_{2}^{T} &=&y_{0}+f\left( \mu _{2}\right) \left( x_{2}-x_{0}\right) \\
&\approx &\underset{y\left( x_{2}\right) }{\underbrace{\left[ y_{0}+f\left(
\mu \left( x_{2}\right) \right) \left( x_{2}-x_{0}\right) \right] }}+\Delta
R_{2}.
\end{eqnarray*}%
The term in square brackets is the exact solution $y\left( x_{2}\right) .$
Therefore, if $\Delta R_{2}$ is suitably small,%
\begin{equation*}
\Delta y_{2}\approx y_{2}^{T}-y_{2},
\end{equation*}%
where $\Delta y_{2}$ denotes the global error in $y_{2}.$ In fact, we
measure the relative global error%
\begin{equation*}
\Delta y_{2}^{rel}\equiv \frac{\Delta y_{2}}{\max (1,\left\vert
y_{2}^{T}\right\vert )}\approx \frac{y_{2}^{T}-y_{2}}{\max (1,\left\vert
y_{2}^{T}\right\vert )}
\end{equation*}%
and we test the inequality%
\begin{equation}
\left\vert \Delta y_{2}^{rel}\right\vert >\left\vert \varepsilon
_{g}-\left\vert \Delta y_{2}^{T}\right\vert \right\vert ,  \label{test >}
\end{equation}%
where $\varepsilon _{g}>0$ is a user-defined tolerance on the relative
global error in the Euler solution, and 
\begin{equation*}
\Delta y_{2}^{T}\equiv \frac{f_{yy}\left( \mu _{2}\right) \left(
x_{2}-x_{0}\right) \Delta _{2}^{2}-2f_{y}\left( \mu _{2}\right) \left(
x_{2}-x_{0}\right) \Delta _{2}}{2\max (1,\left\vert y_{2}^{T}\right\vert )}.
\end{equation*}%
If (\ref{test >}) is true we \textbf{quench} the Euler solution $y_{2}$ by
setting 
\begin{equation*}
y_{2}=y_{0}+f\left( \mu _{2}\right) \left( x_{2}-x_{0}\right) ,
\end{equation*}%
confident that 
\begin{equation*}
\left\vert \Delta y_{2}^{rel}\right\vert \leqslant \max \left( \varepsilon
_{g},\left\vert \Delta y_{2}^{T}\right\vert \right) .
\end{equation*}%
Of course, we prefer that $\left\vert \Delta y_{2}\right\vert $ be bounded
by the imposed tolerance $\varepsilon _{g},$ and it is likely that it will
be if \textbf{strict} local error control has been carried out earlier
(small $\varepsilon _{\rho })$. But even if $\left\vert \Delta
y_{2}^{T}\right\vert >\varepsilon _{g}$ so that $\left\vert \Delta
y_{2}^{T}\right\vert $ is effectively the bound on $\left\vert \Delta
y_{2}^{rel}\right\vert ,$ at least we know the value of $\left\vert \Delta
y_{2}^{T}\right\vert .$ In such a case, if $\left\vert \Delta
y_{2}^{T}\right\vert $ is acceptably small, we would still be satisfied with
the result.

\textbf{NB:} It is imperative to track the growth of $\left\vert \Delta
y_{\cdots }^{T}\right\vert $ and, if $\left\vert \Delta y_{\cdots
}^{T}\right\vert $ becomes too large, we backtrack to a node $x_{k}$ where $%
\left\vert \Delta y_{k}^{T}\right\vert $ is acceptably small, and then
repeat the process from Step 1. That is to say, we treat $\left(
x_{k},y_{k}^{T}\right) $ as a new initial value $\left( x_{0},y_{0}\right) $
and "reboot" the algorithm from this new initial value. This reboot process
is executed \textit{only if it is deemed necessary}.

Obviously, if the test (\ref{test >}) is false, we do not quench the
solution, and simply proceed to the next \ step. An argument can be made
that one should always quench while $\left\vert \Delta y_{\cdots
}^{T}\right\vert $ is acceptably small. We have chosen not to do this in our
numerical examples, in order to make explicit the effects of quenching and
not quenching. Nevertheless, persistent quenching does have its merits.

\item We now have a solution $y_{2}$ that satisfies a tolerance on its
global accuracy. Furthermore, this error control has been implemented in a
stepwise (local) manner. We are ready to compute $y_{3}.$ We return to Step
3, \textbf{update the indices by one} and continue until we eventually reach
the endpoint $x_{N}$.
\end{enumerate}

\section{Numerical Examples}

We consider the following examples for the purpose of numerical experiments:%
\begin{equation*}
\end{equation*}

\begin{center}
\ \ 

\renewcommand{\arraystretch}{1.7}%
\begin{tabular}{|c|c|c|cl}
\multicolumn{4}{l}{\textbf{Table 1:} Details of numerical examples} &  \\ 
\hline
\textbf{\#} & \textbf{ODE} & \textbf{Interval} & $y_{0}$ & 
\multicolumn{1}{|c|}{\textbf{Solution}} \\ \hline\hline
\textbf{1} & \multicolumn{1}{|l|}{$y^{\prime }=y$} & $\left[ 0,5\right] $ & 
\multicolumn{1}{|r}{$2$} & \multicolumn{1}{|l|}{$y\left( x\right) =2e^{x}$}
\\ \hline
\textbf{2} & \multicolumn{1}{|l|}{$y^{\prime }=y^{2}$} & $\left[ -10,-3%
\right] $ & \multicolumn{1}{|r|}{$0.1$} & \multicolumn{1}{|l|}{$y\left(
x\right) =-\frac{1}{x}$} \\ \hline
\textbf{3} & \multicolumn{1}{|l|}{$y^{\prime }=\frac{y}{4}\left( 1-\frac{y}{%
20}\right) $} & $\left[ 0,20\right] $ & \multicolumn{1}{|r|}{$1$} & 
\multicolumn{1}{|l|}{$y\left( x\right) =\frac{20}{1+19e^{-x/4}}$} \\ \hline
\textbf{4} & \multicolumn{1}{|l|}{$y^{\prime }=\frac{1}{y}$} & $\left[ 5,25%
\right] $ & \multicolumn{1}{|r|}{$1$} & \multicolumn{1}{|l|}{$y\left(
x\right) =\sqrt{2x-9}$} \\ \hline
\textbf{5} & \multicolumn{1}{|l|}{$y^{\prime }=\cos y$} & $\left[ a,b\right] 
$ & \multicolumn{1}{|r|}{$-1$} & \multicolumn{1}{|l|}{$x=\ln \left( \sec
y+\tan y\right) $} \\ \hline
\textbf{6} & \multicolumn{1}{|l|}{$y^{\prime }=-y$} & $\left[ 0,10\right] $
& \multicolumn{1}{|r|}{$1$} & \multicolumn{1}{|l|}{$y\left( x\right) =e^{-x}$%
} \\ \hline
\end{tabular}%
\renewcommand{\arraystretch}{1}%
\begin{equation*}
\end{equation*}
\end{center}

where, in \#5,%
\begin{eqnarray*}
a &=&-1.2261911708835170708130609674719 \\
b &=&1.2261911708835170708130609674719.
\end{eqnarray*}%
Some of these problems are weakly stiff, and the solution in \#5 is given
implicitly.

Our results are summarized in Tables $2-8$. In these tables we have used the
notation%
\begin{equation*}
\end{equation*}%
\begin{equation}
\text{ \ \ \ }%
\begin{array}{ccl}
\max \Delta & \text{ } & 
\begin{array}{l}
\text{largest value of }\left\vert \Delta y_{i}^{rel}\right\vert \text{ on
the relevant interval.} \\ 
\end{array}
\\ 
N &  & 
\begin{array}{l}
\text{Number of nodes used.} \\ 
\end{array}
\\ 
\text{Q} &  & 
\begin{array}{l}
\text{Number of nodes where quenching was needed.} \\ 
\end{array}%
\end{array}
\notag
\end{equation}%
\begin{equation}
\text{ \ \ \ }%
\begin{array}{ccl}
\text{P} &  & 
\begin{array}{l}
\text{Number of nodes where local error control in }\mu _{i}\text{ } \\ 
\text{was primary.}%
\end{array}
\\ 
\text{S} &  & 
\begin{array}{l}
\text{Number of nodes where local error control in }\mu _{i}\text{ } \\ 
\text{was secondary.}%
\end{array}
\\ 
\bigstar &  & 
\begin{array}{l}
\text{Number of steps where stepsize control was due to stability } \\ 
\text{considerations (see the later section \textbf{Stiffness and Stability}%
).}%
\end{array}
\\ 
h_{2} &  & 
\begin{array}{l}
\text{size of }h_{2}=x_{2}-x_{1},\text{rounded to first significant digit.}
\\ 
\multicolumn{1}{c}{}%
\end{array}
\\ 
\max \Delta _{E} &  & 
\begin{array}{l}
\text{largest value of }\left\vert \Delta y_{i}^{rel}\right\vert \text{
using the }N\text{ nodes in} \\ 
\text{RK1 without any error control at all.}%
\end{array}%
\end{array}
\notag
\end{equation}%
\begin{equation*}
\end{equation*}

\begin{tabular}{|c|c|c|c|c|c|c|cl}
\multicolumn{8}{l}{\textbf{Table 2:} Numerical results for $\varepsilon
_{g}=10^{-2},\varepsilon _{\rho }=10^{-4}$} &  \\ \hline
\textbf{\#} & $\max \Delta $ & $N$ & \textbf{Q} & \textbf{P} & \textbf{S} & $%
\bigstar $ & $h_{2}$ & \multicolumn{1}{|c|}{$\max \Delta _{E}$} \\ 
\hline\hline
\textbf{1} & $9.3\times 10^{-3}$ & $71$ & $16$ & $0$ & $0$ & $1$ & $%
1.4\times 10^{-3}$ & \multicolumn{1}{|c|}{$2.0\times 10^{-1}$} \\ \hline
\textbf{2} & $8.8\times 10^{-3}$ & $91$ & $1$ & $0$ & $0$ & $1$ & $1.4\times
10^{-3}$ & \multicolumn{1}{|c|}{$1.3\times 10^{-2}$} \\ \hline
\textbf{3} & $9.4\times 10^{-3}$ & $221$ & $1$ & $0$ & $0$ & $1$ & $%
1.4\times 10^{-3}$ & \multicolumn{1}{|c|}{$1.2\times 10^{-2}$} \\ \hline
\textbf{4} & $6.3\times 10^{-3}$ & $221$ & $0$ & $0$ & $0$ & $1$ & $%
1.4\times 10^{-3}$ & \multicolumn{1}{|c|}{$6.3\times 10^{-3}$} \\ \hline
\textbf{5} & $8.9\times 10^{-3}$ & $46$ & $2$ & $0$ & $0$ & $1$ & $1.4\times
10^{-3}$ & \multicolumn{1}{|c|}{$2.1\times 10^{-2}$} \\ \hline
\textbf{6} & $9.2\times 10^{-3}$ & $121$ & $2$ & $0$ & $0$ & $1$ & $%
1.4\times 10^{-3}$ & \multicolumn{1}{|c|}{$1.5\times 10^{-2}$} \\ \hline
\end{tabular}%
\begin{equation*}
\end{equation*}

\begin{tabular}{|c|c|c|c|c|c|c|cl}
\multicolumn{8}{l}{\textbf{Table 3:} Numerical results for $\varepsilon
_{g}=10^{-4},\varepsilon _{\rho }=10^{-6}$} &  \\ \hline
\textbf{\#} & $\max \Delta $ & $N$ & \textbf{Q} & \textbf{P} & \textbf{S} & $%
\bigstar $ & $h_{2}$ & \multicolumn{1}{|c|}{$\max \Delta _{E}$} \\ 
\hline\hline
\textbf{1} & $8.0\times 10^{-5}$ & $88$ & $75$ & $63$ & $0$ & $1$ & $%
1.4\times 10^{-3}$ & \multicolumn{1}{|c|}{$1.6\times 10^{-1}$} \\ \hline
\textbf{2} & $9.2\times 10^{-5}$ & $91$ & $32$ & $0$ & $0$ & $1$ & $%
1.4\times 10^{-3}$ & \multicolumn{1}{|c|}{$1.3\times 10^{-2}$} \\ \hline
\textbf{3} & $9.4\times 10^{-5}$ & $221$ & $115$ & $0$ & $0$ & $1$ & $%
1.4\times 10^{-3}$ & \multicolumn{1}{|c|}{$1.2\times 10^{-2}$} \\ \hline
\textbf{4} & $9.5\times 10^{-5}$ & $223$ & $64$ & $11$ & $0$ & $1$ & $%
1.4\times 10^{-3}$ & \multicolumn{1}{|c|}{$5.6\times 10^{-3}$} \\ \hline
\textbf{5} & $7.9\times 10^{-5}$ & $63$ & $47$ & $30$ & $8$ & $1$ & $%
1.4\times 10^{-3}$ & \multicolumn{1}{|c|}{$1.5\times 10^{-2}$} \\ \hline
\textbf{6} & $9.4\times 10^{-5}$ & $121$ & $53$ & $63$ & $0$ & $1$ & $%
1.4\times 10^{-3}$ & \multicolumn{1}{|c|}{$1.5\times 10^{-2}$} \\ \hline
\end{tabular}%
\begin{equation*}
\end{equation*}

\begin{tabular}{|c|c|c|c|c|c|c|cl}
\multicolumn{8}{l}{\textbf{Table 4:} Numerical results for $\varepsilon
_{g}=10^{-6},\varepsilon _{\rho }=10^{-8}$} &  \\ \hline
\textbf{\#} & $\max \Delta $ & $N$ & \textbf{Q} & \textbf{P} & \textbf{S} & $%
\bigstar $ & $h_{2}$ & \multicolumn{1}{|c|}{$\max \Delta _{E}$} \\ 
\hline\hline
\textbf{1} & $7.2\times 10^{-7}$ & $327$ & $324$ & $310$ & $0$ & $1$ & $%
1.2\times 10^{-3}$ & \multicolumn{1}{|c|}{$3.9\times 10^{-2}$} \\ \hline
\textbf{2} & $7.2\times 10^{-7}$ & $108$ & $91$ & $0$ & $41$ & $1$ & $%
1.4\times 10^{-3}$ & \multicolumn{1}{|c|}{$1.0\times 10^{-2}$} \\ \hline
\textbf{3} & $9.1\times 10^{-7}$ & $221$ & $211$ & $0$ & $0$ & $1$ & $%
1.4\times 10^{-3}$ & \multicolumn{1}{|c|}{$1.2\times 10^{-2}$} \\ \hline
\textbf{4} & $5.4\times 10^{-7}$ & $311$ & $308$ & $142$ & $0$ & $1$ & $%
1.0\times 10^{-3}$ & \multicolumn{1}{|c|}{$2.1\times 10^{-3}$} \\ \hline
\textbf{5} & $8.9\times 10^{-7}$ & $214$ & $209$ & $158$ & $38$ & $1$ & $%
1.4\times 10^{-3}$ & \multicolumn{1}{|c|}{$3.9\times 10^{-3}$} \\ \hline
\textbf{6} & $9.2\times 10^{-7}$ & $276$ & $251$ & $40$ & $219$ & $1$ & $%
1.2\times 10^{-3}$ & \multicolumn{1}{|c|}{$4.0\times 10^{-3}$} \\ \hline
\end{tabular}%
\begin{equation*}
\end{equation*}

\begin{tabular}{|c|c|c|c|c|c|c|cl}
\multicolumn{8}{l}{\textbf{Table 5:} Numerical results for $\varepsilon
_{g}=10^{-8},\varepsilon _{\rho }=10^{-10}$} &  \\ \hline
\textbf{\#} & $\max \Delta $ & $N$ & \textbf{Q} & \textbf{P} & \textbf{S} & $%
\bigstar $ & $h_{2}$ & \multicolumn{1}{|c|}{$\max \Delta _{E}$} \\ 
\hline\hline
\textbf{1} & $8.3\times 10^{-15}$ & $1474$ & $1472$ & $1459$ & $0$ & $1$ & $%
2.7\times 10^{-4}$ & \multicolumn{1}{|c|}{$8.6\times 10^{-3}$} \\ \hline
\textbf{2} & $7.7\times 10^{-9}$ & $330$ & $323$ & $86$ & $219$ & $0$ & $%
1.0\times 10^{-3}$ & \multicolumn{1}{|c|}{$3.4\times 10^{-3}$} \\ \hline
\textbf{3} & $7.9\times 10^{-9}$ & $695$ & $692$ & $663$ & $0$ & $0$ & $%
5.4\times 10^{-4}$ & \multicolumn{1}{|c|}{$4.8\times 10^{-3}$} \\ \hline
\textbf{4} & $1.1\times 10^{-15}$ & $951$ & $949$ & $890$ & $0$ & $0$ & $%
3.3\times 10^{-4}$ & \multicolumn{1}{|c|}{$6.5\times 10^{-4}$} \\ \hline
\textbf{5} & $8.9\times 10^{-10}$ & $949$ & $946$ & $754$ & $176$ & $1$ & $%
3.5\times 10^{-4}$ & \multicolumn{1}{|c|}{$9.0\times 10^{-4}$} \\ \hline
\textbf{6} & $9.5\times 10^{-9}$ & $1233$ & $1223$ & $206$ & $1012$ & $1$ & $%
2.6\times 10^{-4}$ & \multicolumn{1}{|c|}{$8.9\times 10^{-4}$} \\ \hline
\end{tabular}%
\begin{equation*}
\end{equation*}

\begin{tabular}{|c|c|c|c|c|c|c|cl}
\multicolumn{8}{l}{\textbf{Table 6:} Numerical results for $\varepsilon
_{g}=10^{-10},\varepsilon _{\rho }=10^{-12}$} &  \\ \hline
\textbf{\#} & $\max \Delta $ & $N$ & \textbf{Q} & \textbf{P} & \textbf{S} & $%
\bigstar $ & $h_{2}$ & \multicolumn{1}{|c|}{$\max \Delta _{E}$} \\ 
\hline\hline
\textbf{1} & $1.5\times 10^{-14}$ & $6776$ & $6774$ & $6766$ & $0$ & $0$ & $%
1.1\times 10^{-4}$ & \multicolumn{1}{|c|}{$1.9\times 10^{-3}$} \\ \hline
\textbf{2} & $4.7\times 10^{-11}$ & $1460$ & $1457$ & $425$ & $1011$ & $0$ & 
$2.2\times 10^{-4}$ & \multicolumn{1}{|c|}{$7.6\times 10^{-4}$} \\ \hline
\textbf{3} & $1.8\times 10^{-15}$ & $3139$ & $3137$ & $3115$ & $0$ & $0$ & $%
1.2\times 10^{-4}$ & \multicolumn{1}{|c|}{$1.1\times 10^{-3}$} \\ \hline
\textbf{4} & $1.2\times 10^{-15}$ & $4361$ & $4359$ & $4278$ & $0$ & $0$ & $%
6.1\times 10^{-5}$ & \multicolumn{1}{|c|}{$1.5\times 10^{-4}$} \\ \hline
\textbf{5} & $8.6\times 10^{-11}$ & $4354$ & $4350$ & $3518$ & $820$ & $0$ & 
$1.4\times 10^{-4}$ & \multicolumn{1}{|c|}{$2.0\times 10^{-4}$} \\ \hline
\textbf{6} & $2.2\times 10^{-15}$ & $5676$ & $5674$ & $975$ & $4691$ & $0$ & 
$1.0\times 10^{-4}$ & \multicolumn{1}{|c|}{$1.9\times 10^{-4}$} \\ \hline
\end{tabular}%
\begin{equation*}
\end{equation*}

\begin{tabular}{|c|c|c|c|c|c|c|cl}
\multicolumn{8}{l}{\textbf{Table 7:} Numerical results for $\varepsilon
_{g}=10^{-2},\varepsilon _{\rho }=10^{-3}$} &  \\ \hline
\textbf{\#} & $\max \Delta $ & $N$ & \textbf{Q} & \textbf{P} & \textbf{S} & $%
\bigstar $ & $h_{2}$ & \multicolumn{1}{|c|}{$\max \Delta _{E}$} \\ 
\hline\hline
\textbf{1} & $9.3\times 10^{-3}$ & $71$ & $16$ & $0$ & $0$ & $1$ & $%
1.4\times 10^{-3}$ & \multicolumn{1}{|c|}{$2.0\times 10^{-1}$} \\ \hline
\textbf{2} & $8.8\times 10^{-3}$ & $91$ & $1$ & $0$ & $0$ & $1$ & $1.4\times
10^{-3}$ & \multicolumn{1}{|c|}{$1.3\times 10^{-2}$} \\ \hline
\textbf{3} & $9.4\times 10^{-3}$ & $221$ & $1$ & $0$ & $0$ & $1$ & $%
1.4\times 10^{-3}$ & \multicolumn{1}{|c|}{$1.2\times 10^{-2}$} \\ \hline
\textbf{4} & $6.3\times 10^{-3}$ & $221$ & $0$ & $0$ & $0$ & $1$ & $%
1.4\times 10^{-3}$ & \multicolumn{1}{|c|}{$6.3\times 10^{-3}$} \\ \hline
\textbf{5} & $8.9\times 10^{-3}$ & $46$ & $2$ & $0$ & $0$ & $1$ & $1.4\times
10^{-3}$ & \multicolumn{1}{|c|}{$2.1\times 10^{-2}$} \\ \hline
\textbf{6} & $9.2\times 10^{-3}$ & $121$ & $2$ & $0$ & $0$ & $1$ & $%
1.4\times 10^{-3}$ & \multicolumn{1}{|c|}{$1.5\times 10^{-2}$} \\ \hline
\end{tabular}%
\begin{equation*}
\end{equation*}

\begin{tabular}{|c|c|c|c|c|c|c|cl}
\multicolumn{8}{l}{\textbf{Table 8:} Numerical results for $\varepsilon
_{g}=10^{-6},\varepsilon _{\rho }=10^{-7}$} &  \\ \hline
\textbf{\#} & $\max \Delta $ & $N$ & \textbf{Q} & \textbf{P} & \textbf{S} & $%
\bigstar $ & $h_{2}$ & \multicolumn{1}{|c|}{$\max \Delta _{E}$} \\ 
\hline\hline
\textbf{1} & $9.5\times 10^{-7}$ & $162$ & $159$ & $141$ & $0$ & $1$ & $%
1.2\times 10^{-3}$ & \multicolumn{1}{|c|}{$7.9\times 10^{-2}$} \\ \hline
\textbf{2} & $7.2\times 10^{-7}$ & $92$ & $75$ & $0$ & $7$ & $1$ & $%
1.4\times 10^{-3}$ & \multicolumn{1}{|c|}{$1.2\times 10^{-2}$} \\ \hline
\textbf{3} & $9.1\times 10^{-7}$ & $221$ & $211$ & $0$ & $0$ & $1$ & $%
1.4\times 10^{-3}$ & \multicolumn{1}{|c|}{$1.2\times 10^{-2}$} \\ \hline
\textbf{4} & $9.4\times 10^{-7}$ & $243$ & $240$ & $47$ & $0$ & $1$ & $%
1.4\times 10^{-3}$ & \multicolumn{1}{|c|}{$3.6\times 10^{-3}$} \\ \hline
\textbf{5} & $8.9\times 10^{-7}$ & $110$ & $106$ & $69$ & $17$ & $1$ & $%
1.4\times 10^{-3}$ & \multicolumn{1}{|c|}{$7.9\times 10^{-3}$} \\ \hline
\textbf{6} & $9.4\times 10^{-7}$ & $150$ & $137$ & $14$ & $71$ & $1$ & $%
1.2\times 10^{-3}$ & \multicolumn{1}{|c|}{$8.3\times 10^{-3}$} \\ \hline
\end{tabular}%
\begin{equation*}
\end{equation*}%
In every single case, $\max \Delta $ is less than $\varepsilon _{g}.$ For
any given example, the value of $N$ tends to increase as $\varepsilon _{g}$
is decreased. The number of steps where stepsize control was due to
stability considerations was never more than one for any of the examples.
The value of $1.4\times 10^{-3}$ for $h_{2}$ occurs often; this is explained
in the next section. The number of nodes where quenching was needed tends to
increases as $\varepsilon _{g}$ is decreased. LEC was not needed at all for $%
\varepsilon _{g}=10^{-2},$ but for all other cases LEC was necessary and, in
some cases, at almost every node. Primary LEC occurred more often than
secondary LEC, but both types were definitely present. There were, however,
some cases where there was no secondary LEC even though primary LEC was
needed. We included $\max \Delta _{E}$ purely for interests' sake - clearly $%
\max \Delta $ is less than $\max \Delta _{E}$ for every case, usually by
many orders of magnitude (except \#4 in Table 2, where no quenching or LEC
was required). In Tables $2-6,$ $\varepsilon _{\rho }=\varepsilon _{g}/100,$
and in Tables 7 and 8, $\varepsilon _{\rho }=\varepsilon _{g}/10,$
suggesting that a looser bound on LEC still yields acceptable results (see
comments in the Appendix regarding \textit{initial stepsize}).

In Figures $1-10$ (for convenience, all figures are collected at the end of
the paper, after the Appendix), we show some error curves for two of the
examples (\#3 and \#5) for all five global tolerances considered here. In
all of these figures, we adopt the following conventions: horizontal red
line is $\varepsilon _{g};$ black dots indicate $\left\vert \Delta
y_{i}^{rel}\right\vert ,$ the relative global error estimate in the RK1
solution; green curve indicates $\left\vert \Delta y_{i}^{T}\right\vert ,$
the estimate of the relative global error in $y_{i}^{T}$; magenta circles
indicate the true global error in $y_{i}$ at those nodes where quenching was
not needed; blue circles indicate the true global error in $y_{i}$ at those
nodes where quenching was implemented. Clearly, quenching does not occur
when $\left\vert \Delta y_{i}^{rel}\right\vert $ dips below $\varepsilon
_{g},$ as expected. In all cases, the true global error in $y_{i}$ (either
magenta circles or blue circles) lies below $\varepsilon _{g},$ as desired.
Note that $\left\vert \Delta y_{i}^{T}\right\vert $ lies above the quenched
error, but below $\varepsilon _{g},$ thus serving the purpose of a good
estimator. The exception occurs when $\varepsilon _{g}=10^{-10},$ where $%
\left\vert \Delta y_{i}^{T}\right\vert $ is similar to the quenched error,
and both are of the order of machine precision. This probably exposes the
limit of the algorithm's capabilities, at least on the computational
platform used here \cite{platform}, perhaps indicating that a tolerance this
strict is not practical. We also see that quenching becomes more prevalent
as $\varepsilon _{g}$ decreases, because the RK1 solution cannot satisfy
such strict tolerances.

Although we only show results for two examples, the other examples exhibit
similar behaviour. Of course, details for all the examples are to be found
in the tables.

\section{Reboot}

In Step 9 above, we briefly mentioned the reboot process. In Figure 11, we
show the outcome of this process for example \#3, for two reboot tolerances,
on the interval $\left[ 0,50\right] .$ The reboot tolerance $\varepsilon
_{rb}$ is the maximum allowed magnitude of $\left\vert \Delta
y_{i}^{T}\right\vert ,$ and the reboot process is designed to prevent $%
\left\vert \Delta y_{i}^{T}\right\vert $ from becoming too significant with
respect to $\varepsilon _{g}.$ In the first of these figures, we use $%
\varepsilon _{rb}=10^{-10},$ and in the second we use $\varepsilon
_{rb}=10^{-7}.$ It is clear that no reboots were required in the second
case, but in the first case nine reboots were needed (the reboots occur
whenever $\left\vert \Delta y_{i}^{T}\right\vert $ breaches $\varepsilon
_{rb}=10^{-10})$. Admittedly, the value $\varepsilon _{rb}=10^{-10}$ is
strict, and was chosen for the purpose of demonstration. Ordinarily, we
would probably be satisfied with $\varepsilon _{rb}=\varepsilon _{g}/1000.$
For the record, the case $\varepsilon _{rb}=10^{-7}$ has $%
N=521,Q=125,P=0,S=0 $ and $\bigstar =1,$ and the case $\varepsilon
_{rb}=10^{-10}$ has $N=709,Q=111,P=27$ and $S=0$. Each reboot also required
a single stepsize adjustment based on stability/stiffness (see next
section), giving $\bigstar =10.$

The reboot process is, effectively, a type of quench. The idea is to restart
the integration process at the reboot node $x_{rb}$ with the value of $y$
given by%
\begin{equation*}
y_{rb}=y_{0}+f\left( \mu _{rb}^{V}\right) \left( x_{rb}-x_{0}\right) ,
\end{equation*}%
following which we set new values for the initial point%
\begin{eqnarray*}
x_{0} &=&x_{rb} \\
y_{0} &=&y_{rb}
\end{eqnarray*}%
and we proceed from Step 1 described earlier. The reboot node $x_{rb}$ is
the node that precedes the node at which $\left\vert \Delta
y_{i}^{T}\right\vert $ breaches $\varepsilon _{rb}.$ Hence, our reboot
process requires a backtrack of only one node. Note that, since we have
changed the initial point, we need to determine a new $g\left( x,\mu \right)
.$ We necessarily assume that using $\mu _{rb}^{V}$ in $y_{rb}$ will yield
the most accurate value of $y_{rb}$ that we can achieve.

\section{Stiffness and Stability}

We have%
\begin{eqnarray*}
g\left( x,\mu \right) &=&\frac{f\left( y_{0}+f\left( \mu \right) \left(
x-x_{0}\right) \right) -f\left( \mu \right) }{f_{y}\left( \mu \right) \left(
x-x_{0}\right) } \\
&=&\frac{f\left( y_{0}+f\left( \mu \right) \left( x-x_{0}\right) \right)
-f\left( y_{0}+\mu -y_{0}\right) }{f_{y}\left( \mu \right) \left(
x-x_{0}\right) } \\
&\approx &\frac{f\left( y_{0}+f\left( \mu \right) \left( x-x_{0}\right)
\right) -f\left( y_{0}\right) -f_{y}\left( y_{0}\right) \left( \mu
-y_{0}\right) }{\left( \frac{f_{y}\left( \mu \right) }{f\left( \mu \right) }%
\right) f\left( \mu \right) \left( x-x_{0}\right) } \\
&=&\left( \frac{f\left( \mu \right) }{f_{y}\left( \mu \right) }\right)
\left( \frac{f\left( y_{0}+f\left( \mu \right) \left( x-x_{0}\right) \right)
-f\left( y_{0}\right) }{f\left( \mu \right) \left( x-x_{0}\right) }\right) \\
&&+\left( \frac{f\left( \mu \right) }{f_{y}\left( \mu \right) }\right)
\left( \frac{f_{y}\left( y_{0}\right) \left( y_{0}-\mu \right) }{f\left( \mu
\right) \left( x-x_{0}\right) }\right) \\
&\approx &\frac{f\left( \mu \right) f_{y}\left( y_{0}\right) }{f_{y}\left(
\mu \right) }+\frac{f_{y}\left( y_{0}\right) \left( y_{0}-\mu \right) }{%
f_{y}\left( \mu \right) \left( x-x_{0}\right) } \\
&=&\frac{f_{y}\left( y_{0}\right) }{f_{y}\left( \mu \right) }\left( f\left(
\mu \right) +\frac{y_{0}-\mu }{x-x_{0}}\right) .
\end{eqnarray*}%
Hence,%
\begin{equation*}
g_{\mu }\left( x,\mu \right) =-\frac{f_{y}\left( y_{0}\right) f_{yy}\left(
\mu \right) }{f_{y}^{2}\left( \mu \right) }\left( f\left( \mu \right) +\frac{%
y_{0}-\mu }{x-x_{0}}\right) +\frac{f_{y}\left( y_{0}\right) }{f_{y}\left(
\mu \right) }\left( f_{y}\left( \mu \right) -\frac{1}{x-x_{0}}\right) .
\end{equation*}%
For $\mu _{1}\approx y_{0}$ we find%
\begin{eqnarray*}
g_{\mu }\left( x_{1},\mu _{1}\right) &\approx &-\frac{f_{y}\left(
y_{0}\right) f_{yy}\left( y_{0}\right) }{f_{y}^{2}\left( y_{0}\right) }%
\left( f\left( y_{0}\right) +\frac{y_{0}-y_{0}}{x_{1}-x_{0}}\right) \\
&&+\frac{f_{y}\left( y_{0}\right) }{f_{y}\left( y_{0}\right) }\left(
f_{y}\left( y_{0}\right) -\frac{1}{x_{1}-x_{0}}\right) \\
&=&-\frac{f_{yy}\left( y_{0}\right) f\left( y_{0}\right) }{f_{y}\left(
y_{0}\right) }+f_{y}\left( y_{0}\right) -\frac{1}{x_{1}-x_{0}}
\end{eqnarray*}%
and if $x_{1}\approx x_{0}$ and $y^{\prime }=f\left( y\right) $ is nonstiff
or weakly stiff, we would likely find that $\frac{-1}{x_{1}-x_{0}}$ is the
dominant term. So,%
\begin{equation*}
g_{\mu }\left( x_{1},\mu _{1}\right) \approx -\left( \frac{1}{x_{1}-x_{0}}%
\right) =-1000
\end{equation*}%
for the values used in this study. We note that if $f_{y}\left( y_{0}\right) 
$ is small and positive, the first term on the RHS could contribute to
stiffness in $\mu ^{\prime }=g\left( x,\mu \right) ;$ if $f_{y}\left(
y_{0}\right) $ is negative the second term contributes to stiffness in $\mu
^{\prime }=g\left( x,\mu \right) ;$ if $f_{y}\left( y_{0}\right) $ is
negative and has suitably small magnitude, the first term could cancel
somewhat with the third term, thus reducing the stiffness in $\mu ^{\prime
}=g\left( x,\mu \right) .$

For the examples considered here, we found%
\begin{equation*}
\end{equation*}

\begin{center}
\renewcommand{\arraystretch}{1.15}%
\begin{tabular}{|c|c|c|cl}
\multicolumn{4}{l}{\textbf{Table 9:} $\mu _{1}$ and $g_{\mu }\left(
x_{1},\mu _{1}\right) $} &  \\ \cline{1-4}
\textbf{\#} & $y_{0}$ & $\mu _{1}$ & $g_{\mu }\left( x_{1},\mu _{1}\right) $
& \multicolumn{1}{|l}{} \\ \cline{1-4}
\textbf{1} & \multicolumn{1}{|r|}{$2$} & \multicolumn{1}{|r|}{$2.0010$} & 
\multicolumn{1}{|r}{$-999.00$} & \multicolumn{1}{|l}{} \\ 
\cline{1-4}\cline{4-4}
\textbf{2} & \multicolumn{1}{|r|}{$0.1$} & \multicolumn{1}{|r|}{$0.1000$} & 
\multicolumn{1}{|r}{$-999.85$} & \multicolumn{1}{|l}{} \\ 
\cline{1-4}\cline{4-4}
\textbf{3} & \multicolumn{1}{|r|}{$1$} & \multicolumn{1}{|r|}{$1.0001$} & 
\multicolumn{1}{|r}{$-999.76$} & \multicolumn{1}{|l}{} \\ 
\cline{1-4}\cline{4-4}
\textbf{4} & \multicolumn{1}{|r|}{$1$} & \multicolumn{1}{|r|}{$1.0005$} & 
\multicolumn{1}{|r}{$-999.99$} & \multicolumn{1}{|l}{} \\ 
\cline{1-4}\cline{4-4}
\textbf{5} & \multicolumn{1}{|r|}{$-1$} & \multicolumn{1}{|r|}{$-0.9997$} & 
\multicolumn{1}{|r}{$-998.98$} & \multicolumn{1}{|l}{} \\ 
\cline{1-4}\cline{4-4}
\textbf{6} & \multicolumn{1}{|r|}{$1$} & \multicolumn{1}{|r|}{$0.8885$} & 
\multicolumn{1}{|r}{$-1001.00$} & \multicolumn{1}{|l}{} \\ 
\cline{1-4}\cline{4-4}
\end{tabular}%
\renewcommand{\arraystretch}{1}%
\begin{equation*}
\end{equation*}
\end{center}

consistent with our assumption $\mu _{1}\approx y_{0}$ and our result $%
g_{\mu }\left( x_{1},\mu _{1}\right) \approx -1000.$

For stability, we demand that, when solving the Dahlquist equation%
\begin{equation*}
\mu ^{\prime }=-\lambda \mu ,\text{ }\lambda >0,\text{ }\mu \left( 0\right)
=\mu _{0}=1
\end{equation*}%
using any of the RK methods in DP853, we find a solution at $x_{1}$ that has
the same properties as the true solution, i.e.%
\begin{equation*}
0<e^{-\lambda x_{1}}<1.
\end{equation*}%
In other words, we demand that if 
\begin{equation*}
\mu _{1}=\mu _{0}+h\Phi \left( x_{0},\mu _{0}\right) ,
\end{equation*}%
then%
\begin{equation}
0<\mu _{1}<1.  \label{0 < mu <1}
\end{equation}%
Usually, stability with regard to stiffness concerns only the upper bound in
this expression, but we have observed that the methods in DP853 yield
negative values of $\mu _{1}$ for some values of $h.$ Hence, we have
determined the range of values of $h$ such that both bounds in (\ref{0 < mu
<1}) are satisfied by all three methods in DP853. This has yielded%
\begin{equation*}
-1.3764<-\lambda h<0.
\end{equation*}%
For our situation involving $\mu ^{\prime }=g\left( x,\mu \right) ,$ the
stiffness constant $\lambda $ for the step $\left[ x_{1},x_{2}\right] $ is
estimated as%
\begin{equation*}
\lambda =\left\vert g_{\mu }\left( x_{1},\mu _{1}\right) \right\vert
\end{equation*}%
so that 
\begin{equation}
h_{2}=\frac{1.3764}{\left\vert g_{\mu }\left( x_{1},\mu _{1}\right)
\right\vert }\approx 1.3764\times 10^{-3}  \label{stability h2}
\end{equation}%
is the largest value of $h_{2}=x_{2}-x_{1}$ that we should be willing to
use. This value (rounded to $1.4\times 10^{-3})$ appears numerous times in
the tables of results presented earlier, as the value of $h_{2}$ that was
used. This restriction on the magnitude of the stepsize is performed at
every step \textit{prior} to the numerical solution on that step being
computed. Any stepsize adjustment via LEC that is subsequently required is,
of course, performed \textit{after} the numerical solution on that step has
been found.

Note that we do not need to determine the stability regions of DP853 in the
complex plane, because we are dealing only with univariate problems in this
work. The stability interval $\left[ -1.3764,0\right] $ on the negative real
axis is all that we require. For systems of differential equations, the
stiffness constant would be estimated using the eigenvalues of the Jacobian
of the system, which could be complex numbers. In such a case, the entire
stability region of each method in DP853 in the complex plane would need to
be known.

\section{Future Research}

We regard this paper as the first in a larger project, and so we have
restricted our study to that of scalar, autonomous, nonstiff or weakly stiff
problems. We anticipate that future research will address the following:

\begin{enumerate}
\item Non-autonomous problems $y^{\prime }=f\left( x,y\right) $, which will
require two unknowns, $\xi $ and $\mu ,$ to be determined for the remainder
term. This means that we will need to solve the system%
\begin{equation*}
\left[ 
\begin{array}{c}
\xi ^{\prime } \\ 
\mu ^{\prime }%
\end{array}%
\right] =\left[ 
\begin{array}{c}
g_{1}\left( x,\xi ,\mu \right) \\ 
g_{2}\left( x,\xi ,\mu \right)%
\end{array}%
\right] .
\end{equation*}

\item Non-autonomous systems%
\begin{equation*}
\left[ 
\begin{array}{c}
y_{1}^{\prime } \\ 
\vdots \\ 
y_{m}^{\prime }%
\end{array}%
\right] =\left[ 
\begin{array}{c}
f_{1}\left( x,y_{1},\ldots ,y_{m}\right) \\ 
\vdots \\ 
f_{m}\left( x,y_{1},\ldots ,y_{m}\right)%
\end{array}%
\right] ,
\end{equation*}%
which we anticipate will yield $m+1$ equations for each component, giving $%
m^{2}+m$ equations overall that will need to be solved to find the remainder
term.

\item The use of \textit{tandem} RK methods instead of DP853. For example,
we may use RK2, RK4, and RK7 to generate low-, high- and very high-order
solutions for $\xi $ and $\mu ,$ even if these RK methods are not embedded.
If this proves successful, it means that we are not bound to DP853, and that
might provide flexibility for the algorithm.

\item Strongly stiff systems must be considered, wherein we might have to
utilize the implicit form of RK1 and suitably stable methods in place of
DP853. Indeed, we believe our algorithm has a modular character, so that the
RK methods we have used could be replaced with other methods that may have
properties more suitable to the problem at hand, whilst preserving the
structure of the algorithm as a whole. As a matter of record, we have solved
the stiff problem $y^{\prime }=-1000y$ quite easily using our algorithm,
although very small stepsizes were needed $(h\sim 10^{-3})$ due to the
explicit nature of DP853. Replacing DP853 with \textit{A}-stable methods
might be a useful line of enquiry.

\item The use of low-order RK methods other than Euler's method.

\item Can a stricter estimate of the global error $\left\vert \Delta
y_{i}^{T}\right\vert $ be made, and is it necessary to do so?
\end{enumerate}

\section{Concluding Comments}

We have described a novel algorithm for controlling global error in a
stepwise sense. The algorithm exploits the remainder term of a Taylor
expansion of the solution, and requires the use of the DP853 triple to solve
an auxiliary problem which enables this remainder term to be determined. A
quenching process then allows the solution generated by Euler's method to be
controlled. We have achieved tolerances on the relative global error as
strict as $10^{-10}.$ Admittedly, we have restricted the calculations in
this paper to scalar, autonomous, nonstiff or weakly stiff problems.
Naturally, we hope that future research will expand the scope of our
algorithm.

\section{Appendix}

\subsection{Newton's Method for $\protect\mu _{1}$}

We must solve%
\begin{equation*}
y_{\mu }-y_{0}-f\left( \mu _{1}\right) \left( x_{\mu }-x_{0}\right) =0
\end{equation*}%
for $\mu _{1}.$ It is not enough to merely solve this equation in the form
presented here. Recall from (\ref{Taylor}) that we actually have 
\begin{equation*}
y\left( x_{\mu }\right) =y\left( x_{0}\right) +y^{\prime }\left( \xi
_{1},\mu _{1}\right) \left( x_{\mu }-x_{0}\right) ,
\end{equation*}%
where%
\begin{equation*}
\mu _{1}\equiv y\left( \xi _{1}\right) .
\end{equation*}%
Strictly speaking, then, we must solve%
\begin{equation*}
F_{\xi }\left( \xi _{1}\right) \equiv y_{\mu }-y_{0}-y^{\prime }\left( \xi
_{1},y\left( \xi _{1}\right) \right) \left( x_{\mu }-x_{0}\right) =0
\end{equation*}%
to find $\xi _{1}.$ We use Newton's Method in the following manner:%
\begin{eqnarray*}
&&\text{Use RK7 on }\left[ x_{0},\xi _{1}^{j}\right] \text{ to estimate }%
y\left( \xi _{1}^{j}\right) \\
&&\text{Use RK7 on }\left[ x_{0},\xi _{1}^{j}+\delta \xi \right] \text{ to
estimate }y\left( \xi _{1}^{j}+\delta \xi \right) \\
&&F_{\xi }^{\prime }\left( \xi _{1}^{j}\right) =\frac{F_{\xi }\left( \xi
_{1}^{j}+\delta \xi \right) -F_{\xi }\left( \xi _{1}^{j}\right) }{\delta \xi 
} \\
&&\xi _{1}^{j+1}=\xi _{1}^{j}-\frac{F_{\xi }\left( \xi _{1}^{j}\right) }{%
F_{\xi }^{\prime }\left( \xi _{1}^{j}\right) },
\end{eqnarray*}%
where we use $\xi _{1}^{0}=\left( x_{\mu }+x_{0}\right) /2$ as the starting
value for the Newton iteration. We have used $\delta \xi =10^{-5}$ in all
our calculations, and we found that no more than three Newton iterations
were ever needed for all the examples that we have considered. The stopping
criterion we used for the iteration scheme was simply%
\begin{equation*}
\left\vert F_{\xi }\left( \xi _{1}^{j}\right) \right\vert <10^{-14}.
\end{equation*}%
Once we have found an acceptable value for $\xi _{1},$ we use RK7 on $\left[
x_{0},\xi _{1}\right] $ to estimate $y\left( \xi _{1}\right) ,$ and we take
this estimate to be $\mu _{1}.$ In the implementations of RK7 here, we used
five steps on the relevant intervals. We note that if the problem is
strongly stiff, we would probably prefer to use an \textit{A}-stable RK
method instead of RK7.

\subsection{Suggested local tolerance and initial stepsize}

\subsubsection{Local tolerance}

Obtain a "quick and dirty" solution for $\mu ^{\prime }=g\left( x,\mu
\right) $ using an RK method of moderate order, such as RK4 (which has only
four internal stages, and so is reasonably efficient for this purpose).
Denote this solution $\left( x_{j},\mu _{j}\right) .$ Assume that the
resulting relative global error in $y_{j}=y_{0}+f\left( \mu _{j}\right)
\left( x_{j}-x_{0}\right) $ is given by%
\begin{equation*}
\Delta _{j}^{rel}=\frac{f_{y}\left( \mu _{j}\right) \left(
x_{j}-x_{0}\right) \delta \mu _{j}}{\max\nolimits_{j}\left( 1,\left\vert
y_{0}+R\left( x_{j},\mu _{j}\right) \right\vert \right) }=\frac{f_{y}\left(
\mu _{j}\right) \left( x_{j}-x_{0}\right) \delta \mu _{j}}{%
\max\nolimits_{j}\left( 1,\left\vert y_{0}+f\left( \mu _{j}\right) \left(
x_{j}-x_{0}\right) \right\vert \right) }
\end{equation*}%
where $\delta \mu _{j}$ denotes the global error in $\mu _{j}.$ We demand
that this error be equal in magnitude, at most, to the specified global
tolerance $\varepsilon _{g}.$ This gives%
\begin{align*}
& \left\vert \frac{f_{y}\left( \mu _{j}\right) \left( x_{j}-x_{0}\right)
\delta \mu _{j}}{\max\nolimits_{j}\left( 1,\left\vert y_{0}+f\left( \mu
_{j}\right) \left( x_{j}-x_{0}\right) \right\vert \right) }\right\vert
\leqslant \varepsilon _{g} \\
\Rightarrow & \left\vert \delta \mu _{j}\right\vert \leqslant \frac{%
\varepsilon _{g}\max\nolimits_{j}\left( 1,\left\vert y_{0}+f\left( \mu
_{j}\right) \left( x_{j}-x_{0}\right) \right\vert \right) }{\left\vert
f_{y}\left( \mu _{j}\right) \left( x_{j}-x_{0}\right) \right\vert }.
\end{align*}%
Choose 
\begin{equation*}
\delta \mu _{m}\equiv \max\nolimits_{j}\left\vert \frac{\varepsilon
_{g}\max\nolimits_{j}\left( 1,\left\vert y_{0}+f\left( \mu _{j}\right)
\left( x_{j}-x_{0}\right) \right\vert \right) }{\left\vert f_{y}\left( \mu
_{j}\right) \left( x_{j}-x_{0}\right) \right\vert }\right\vert .
\end{equation*}%
Assume 
\begin{equation*}
\delta \mu _{m}=\left( x_{N}-x_{0}\right) \varepsilon _{\rho }
\end{equation*}%
and so, using (\ref{eps rho <= Delta/interval}),%
\begin{equation*}
\varepsilon _{\rho }=\frac{\delta \mu _{m}}{\left( x_{N}-x_{0}\right) }
\end{equation*}%
is an estimate for the value of $\varepsilon _{\rho }$ consistent with the
specified global tolerance $\varepsilon _{g}.$ We also consider a
user-defined default value for $\varepsilon _{\rho },$ denoted $\varepsilon
_{\rho }^{D},$ and we choose the local tolerance via%
\begin{equation*}
\varepsilon _{\rho }=\min \left( \varepsilon _{\rho }^{D},\frac{\delta \mu
_{m}}{\left( x_{N}-x_{0}\right) }\right) .
\end{equation*}%
In this study we have generally used $\varepsilon _{\rho }^{D}=\varepsilon
_{g}/100,$ but we have done some calculations using $\varepsilon _{\rho
}^{D}=\varepsilon _{g}/10$ (see Tables 7 and 8).

Lastly, we note that if the problem is strongly stiff, we would probably
prefer to use an \textit{A}-stable RK method instead of RK4.

\subsubsection{Initial stepsize}

Assuming that DP3 is formally equivalent to a Taylor expansion up to
third-order, we treat the fourth-order term%
\begin{equation*}
\frac{\mu ^{\left( iv\right) }\left( x_{1},\mu _{1}\right) h_{2}^{4}}{24}
\end{equation*}%
as the numerical error at $x_{2}$ obtained using DP3. The relevant stepsize $%
h_{2}$ is estimated from%
\begin{equation*}
\left\vert \frac{\mu ^{\left( iv\right) }\left( x_{1},\mu _{1}\right)
h_{2}^{4}}{24h_{2}\max (1,\left\vert \mu _{1}\right\vert )}\right\vert
=\left\vert \frac{\mu ^{\left( iv\right) }\left( x_{1},\mu _{1}\right)
h_{2}^{3}}{24\max (1,\left\vert \mu _{1}\right\vert )}\right\vert
=\varepsilon _{\rho }.
\end{equation*}%
Using 
\begin{align*}
\mu ^{\left( iv\right) }=g^{\prime \prime \prime }& =\left( \frac{\partial }{%
\partial x}+\frac{d\mu }{dx}\frac{\partial }{\partial \mu }\right) \left( 
\frac{\partial }{\partial x}+\frac{d\mu }{dx}\frac{\partial }{\partial \mu }%
\right) \left( \frac{\partial }{\partial x}+\frac{d\mu }{dx}\frac{\partial }{%
\partial \mu }\right) g \\
& =\left( \frac{\partial }{\partial x}+g\frac{\partial }{\partial \mu }%
\right) \left( \frac{\partial }{\partial x}+g\frac{\partial }{\partial \mu }%
\right) \left( \frac{\partial }{\partial x}+g\frac{\partial }{\partial \mu }%
\right) g
\end{align*}%
gives%
\begin{equation}
h_{2}=\left\vert \frac{24\max (1,\left\vert \mu _{1}\right\vert )}{g^{\prime
\prime \prime }\left( x_{1},\mu _{1}\right) }\right\vert ^{\frac{1}{3}}.
\label{h2 g'''}
\end{equation}%
Naturally, we also specify a default value $h_{2}^{D},$ and we invoke the
stability consideration (\ref{stability h2}), so that%
\begin{equation*}
h_{2}=\min \left( h_{2}^{D},\frac{1.3764}{\left\vert g_{\mu }\left(
x_{1},\mu _{1}\right) \right\vert },\left\vert \frac{24\max (1,\left\vert
\mu _{1}\right\vert )}{g^{\prime \prime \prime }\left( x_{1},\mu _{1}\right) 
}\right\vert ^{\frac{1}{3}}\right) .
\end{equation*}%
In our calculations we used $h_{2}^{D}=0.1,$ and found that the initial
stepsize was always set either by (\ref{stability h2}) or (\ref{h2 g'''}).

\begin{equation*}
\end{equation*}

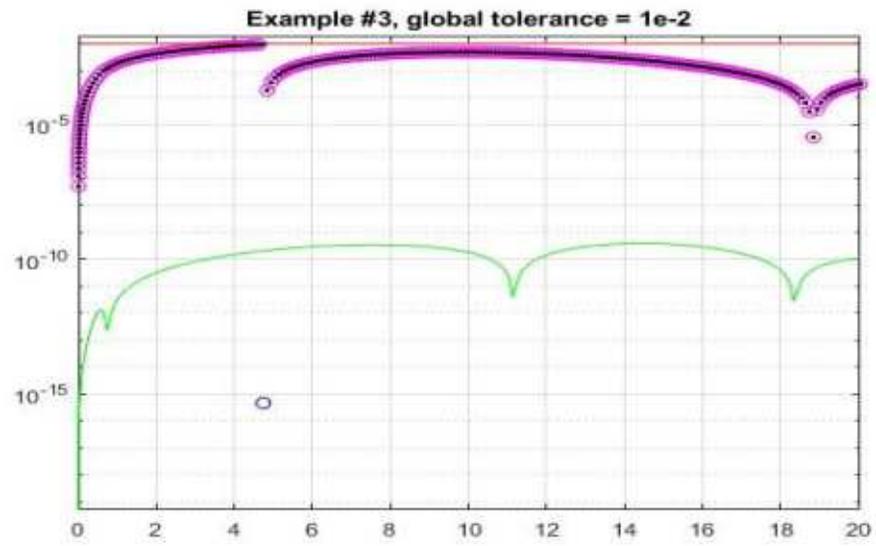
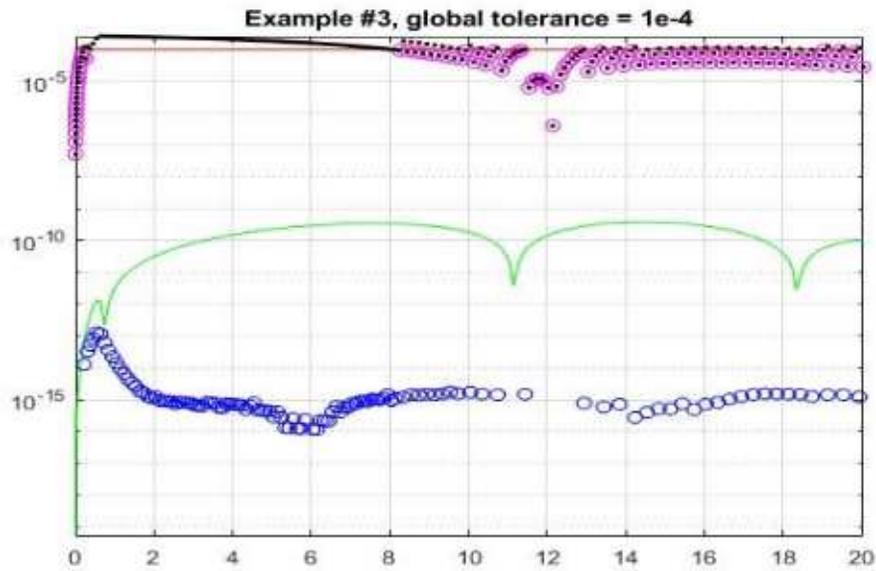
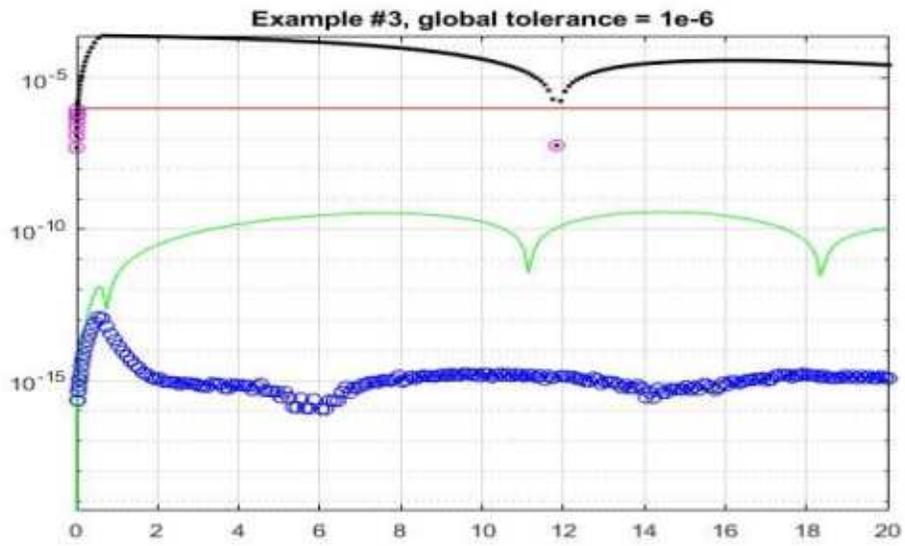
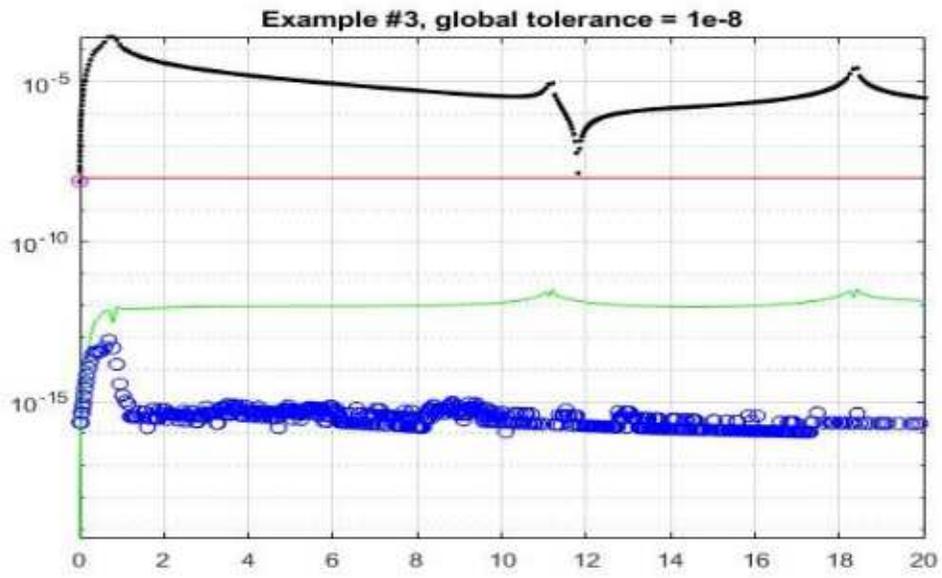
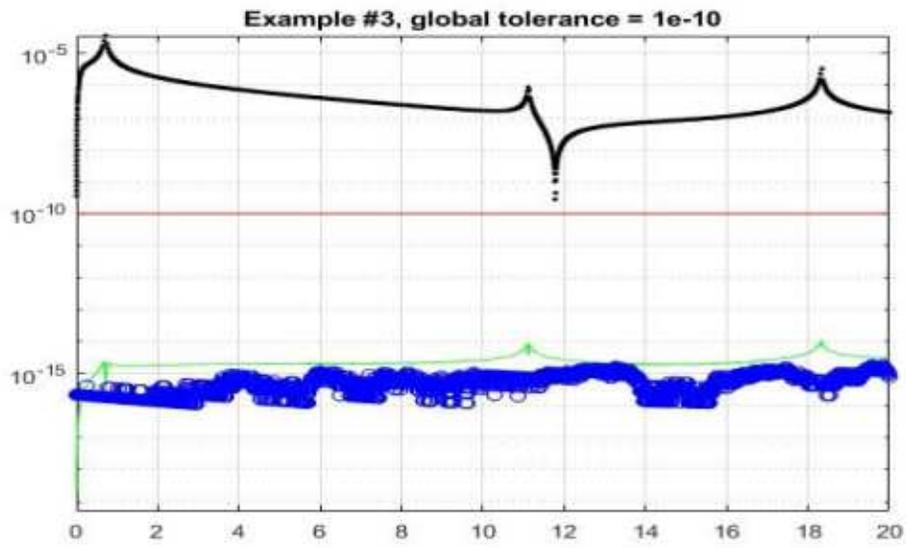
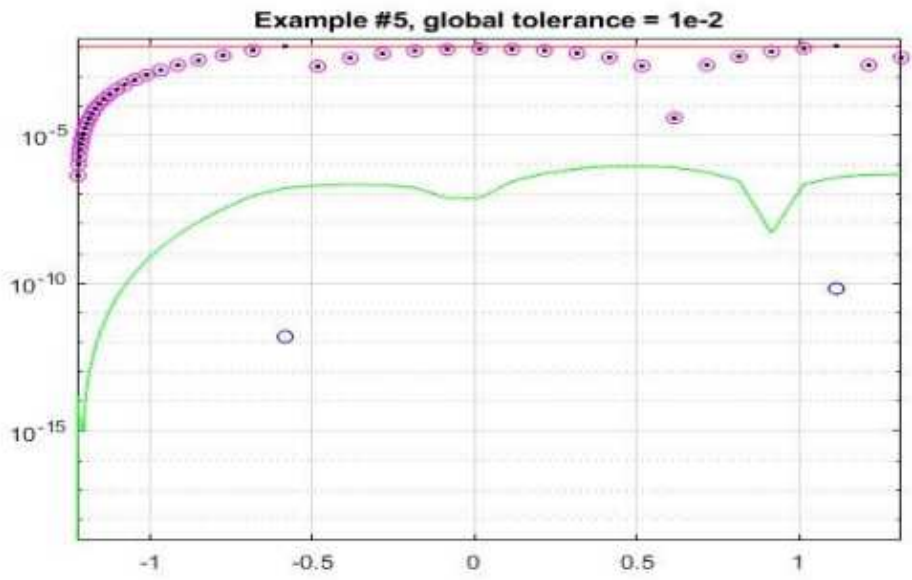
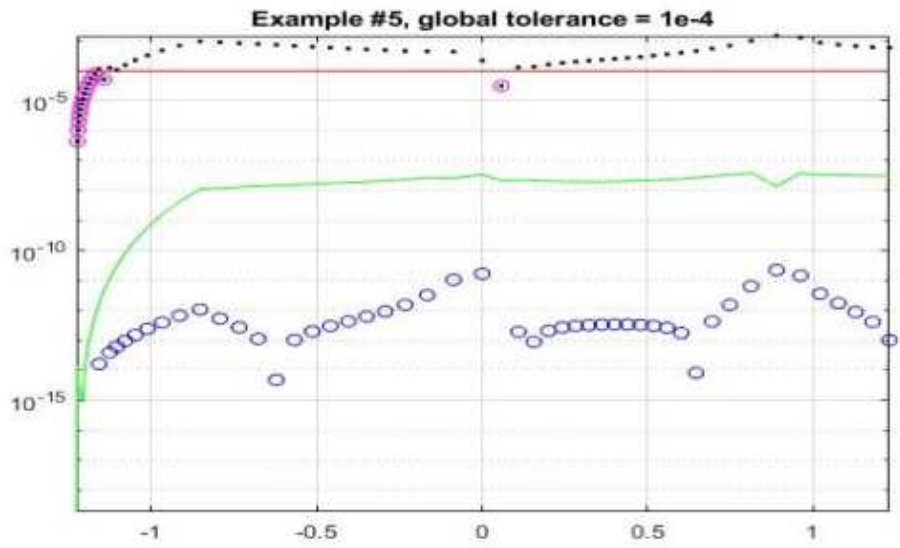
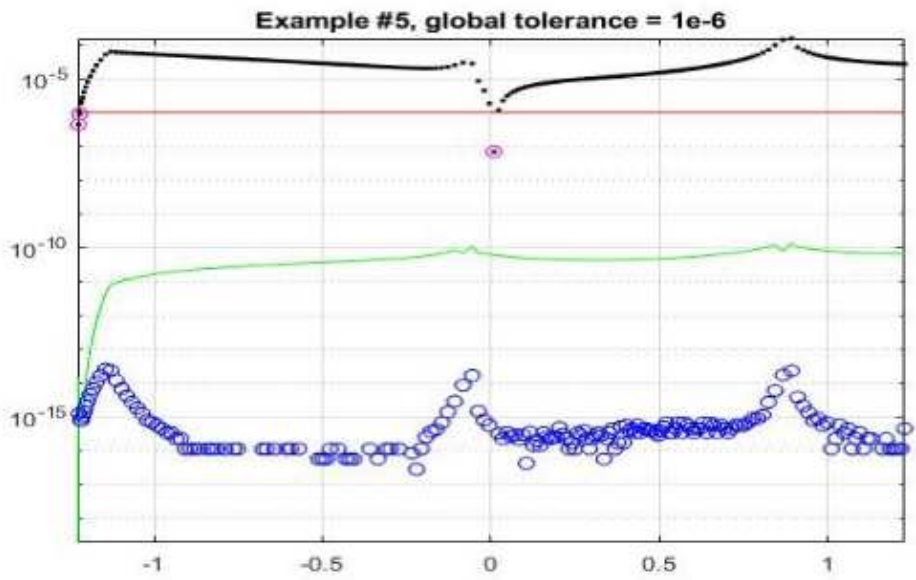
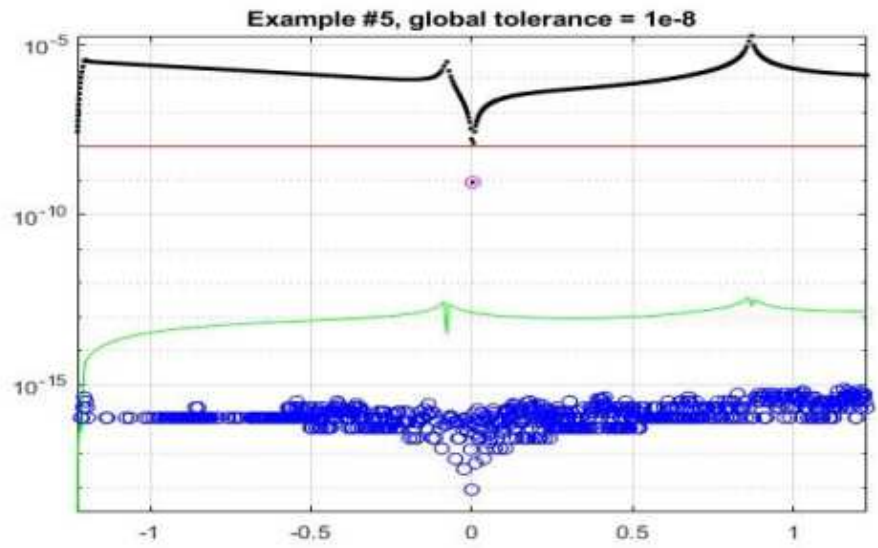
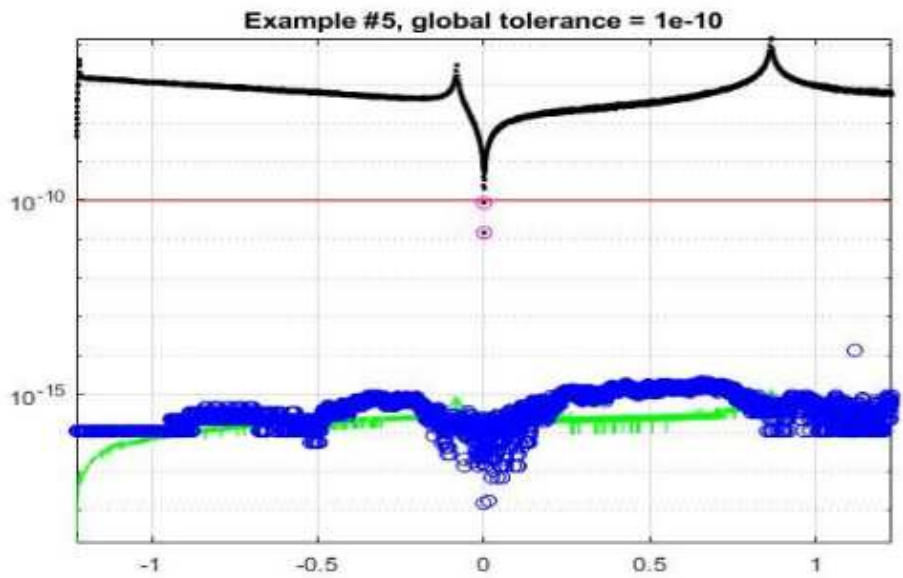
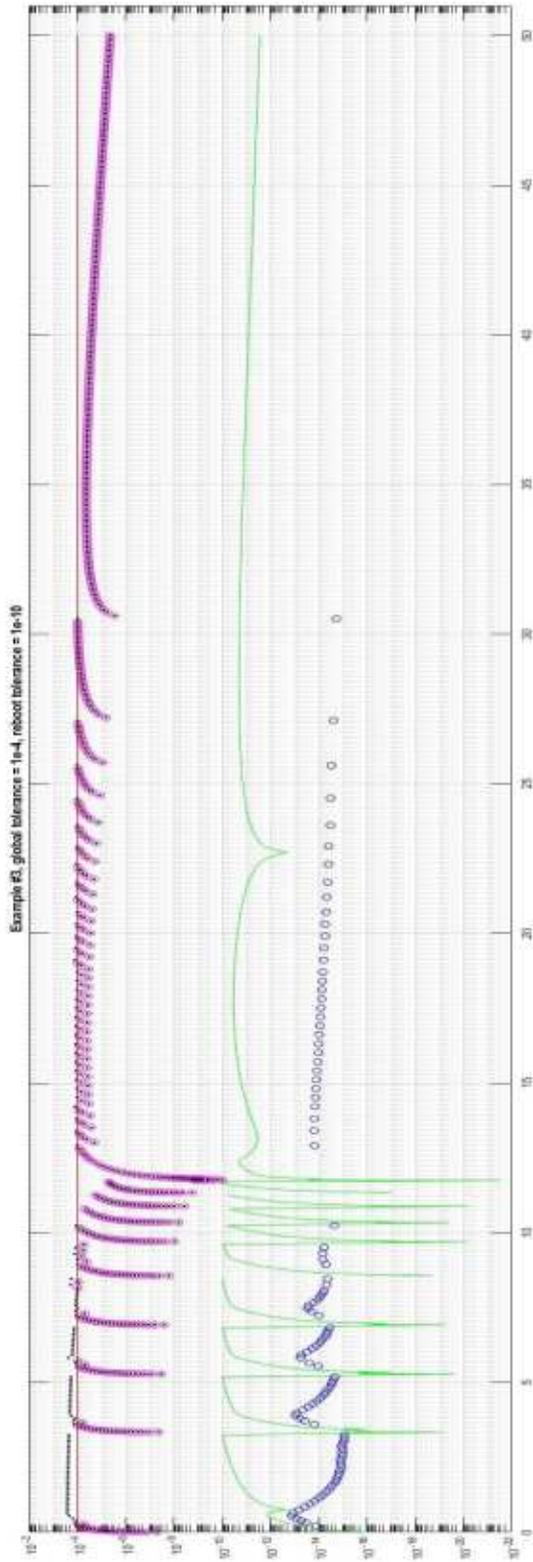
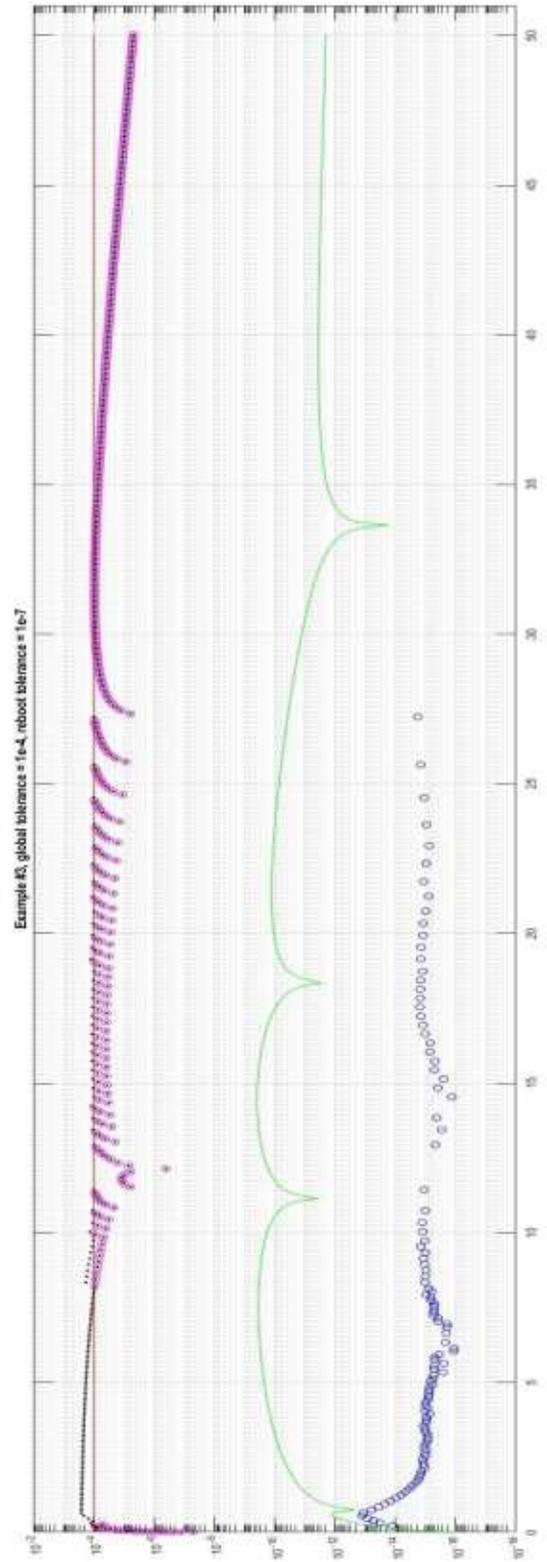
\end{document}